\input  amstex
\input amsppt.sty
\magnification1200
\vsize=23.5truecm
\hsize=16.5truecm
\NoBlackBoxes
\def\d{d\!@!@!@!@!@!{}^{@!@!\text{\rm--}}\!}

\def\supp{\operatorname{supp}}

\def\crp{\overline{\Bbb R}_+}

\def\crpm{\overline{\Bbb R}_\pm}

\def\rnp{{\Bbb R}^n_+}

\def\rnpm{\Bbb R^n_\pm}
\def\crnp{\overline{\Bbb R}^n_+}
\def\crnm{\overline{\Bbb R}^n_-}
\def\crnpm{\overline{\Bbb R}^n_\pm}
\def\comega{\overline\Omega }

\def\ang#1{\langle {#1} \rangle}

\def\Op{\operatorname{Op}}
\def\Pfrac{\tsize\frac1{\raise 1pt\hbox{$\scriptstyle p$}}}
\def\pfrac{\frac1{\raise 1pt\hbox{$\scriptscriptstyle p$}}}
\def\Pfracc#1{\tsize\frac{#1}{\raise 1pt\hbox{$\scriptstyle p$}}}
\def\pfracc#1{\frac{#1}{\raise 1pt\hbox{$\scriptscriptstyle p$}}}

\def\simto{\overset\sim\to\rightarrow}

\def\Zfrac{\tsize\frac1{\raise 1pt\hbox{$\scriptstyle z$}}}
\def\zfrac{\frac1{\raise 1pt\hbox{$\scriptscriptstyle z$}}}

\def\rp{ \Bbb R_+}

\def\rmi{ \Bbb R_-}

\def\R{\Bbb R}

\def\ol{\overline}
\def\SD{\Cal S}
\def\E{\Cal E}
\def\F{\Cal F}
\def\D{\Cal D}

\document
\topmatter
\title
Integration by parts and  Pohozaev identities for space-dependent fractional-order operators
\endtitle
\author Gerd Grubb \endauthor
\affil
{Department of Mathematical Sciences, Copenhagen University,
Universitetsparken 5, DK-2100 Copenhagen, Denmark.
E-mail {\tt grubb\@math.ku.dk}}\endaffil

\abstract
Consider a classical elliptic  pseudodifferential operator $P$ on ${\Bbb R}^n$ of order $2a$
($0<a<1)$ with even symbol. For example, $P=A(x,D)^a$ where $A(x,D)$ is a
second-order strongly elliptic differential operator; the fractional Laplacian
$(-\Delta )^a$ is a particular case. For solutions $u$ of the Dirichlet problem on a bounded
smooth subset $\Omega \subset{\Bbb R}^n$, we show an
integration-by-parts formula with a boundary integral involving $(d^{-a}u)|_{\partial\Omega }$, where $d(x)=\operatorname{dist}(x,\partial\Omega )$. This
extends recent results of Ros-Oton, Serra and Valdinoci, to operators
that are $x$-dependent, nonsymmetric, and have lower-order parts. We
also generalize their formula of Pohozaev-type, that can be used to
prove unique continuation properties, and nonexistence of nontrivial
solutions of semilinear problems. An illustration is given with $P=(-\Delta +m^2)^a$. The basic step in our analysis is a
factorization of $P$, $P\sim P^-P^+$, where we set up a calculus for
the generalized pseudodifferential operators $P^\pm$ that come out of the construction. 
\endabstract
\endtopmatter

\rightheadtext {Integration by parts}

\subhead 1. Introduction 
\endsubhead

A prominent example of a fractional-order  pseudodifferential operator
 ($\psi $do) is
the fractional Laplacian $(-\Delta )^a$ on ${\Bbb R}^n$, $0<a<1$; 
$$
(-\Delta )^au=\operatorname{Op}(|\xi |^{2a})u=
\Cal F^{-1}(|\xi |^{2a}\hat u(\xi )),\quad \hat
u(\xi )=\Cal F
u= \int_{{\Bbb R}^n}e^{-ix\cdot \xi }u(x)\, dx.\tag1.1
$$
It is
currently of 
great interest in probability, finance,  mathematical physics and
differential geometry.
It can also be described as a singular integral operator 
$$
(-\Delta )^au(x)=c_{n,a}PV\int_{{\Bbb
R}^n}\frac{u(x)-u(y)}{|x-y|^{n+2a}}\,dy
=c_{n,a}PV\int_{{\Bbb R}^n}\frac{u(x)-u(x+y)}{|y|^{n+2a}}\,dy,\tag1.2
$$
with convolution kernel $c_{n,a}|y|^{-n-2a}=\F^{-1}|\xi |^{2a}$.

Both descriptions allow generalizations. In (1.1), one can replace the
symbol $|\xi |^{2a}$ by a more general nonvanishing function $p_0(x,\xi )\in
C^\infty ({\Bbb R}^n\times({\Bbb R}^n\setminus \{0\}))$ that is
homogeneous in $\xi $ of degree $2a$, and add terms of lower order, to
get a classical $\psi $do symbol $p(x,\xi )$; the operator is then no longer
translation-invariant nor symmetric. Such operators are standard
examples in the pseudodifferential calculus, and their boundary
value problems  on  suitably smooth subsets $\Omega $ of ${\Bbb R}^n$ have been treated in works of Vishik and Eskin, cf.\
e.g.\ \cite{E81},  Duduchava et al.\ \cite{D84,CD01}, and Shargorodsky
\cite{S94}, with results on solvability in limited ranges of Sobolev
spaces. Recently, a new boundary value theory has been presented
in Grubb \cite{G14,G15}, obtaining
regularity estimates of solutions divided by $d^a$
($d(x)=\operatorname{dist}(x,\partial\Omega )$) in full scales of function
spaces with orders $s\to\infty $, for example in
H\"older spaces of arbitrarily high order. 

The pseudodifferential theory is useful in allowing a direct treatment of 
$x$-dependent operators, providing solution operators (or parametrices) that can give
more efficient regularity estimates than the technique of perturbation of constant-coefficient cases.

In (1.2), one can replace the function $c_{n,a}|y|^{-n-2a}$ by other
positive functions $K(y)$ that are homogeneous in $y$ of degree
$-n-2a$ and possibly less smooth. (In the smooth case this coincides
with $\psi $do's with homogeneous $x$-independent symbol.) Such cases and further
generalizations 
have recently been
studied in probability and nonlinear analysis, see e.g.\ Caffarelli and
Silvestre \cite{CS09}, 
Ros-Oton and Serra  \cite{RS14a,RS15b}, and their
references. 
For problems on bounded domains $\Omega $, the integral operator methods allow limited
 smoothness of the integrand and boundary. To our knowledge, they have
 with few exceptions been applied to $x$-independent (translation-invariant) positive selfadjoint operators.

In the generalizations of (1.1) and (1.2), the fact that  $|\xi
|^{2a}$ is {\it even} (takes
the same value at $\xi $ and $-\xi $) is kept as a hypothesis, that
$p_0$ is even in $\xi $, resp.\ that $K$ is even in $y$.

The methods used in the pseudodifferential theory are complex, and
differ radically
from the real methods currently used for the singular integral
formulations.

There is a large number of preceding studies of boundary problems for
$(-\Delta )^a$ and its generalizations; let us mention e.g.\
\cite{BG59,L72,HJ96,K97,CS98,J02,S07,MN14,SV14}, \cite{BBS15,FG15,FKV15,BSV15}.

\medskip

A useful tool in solvability studies for linear and nonlinear
partial differential equations on subsets $\Omega \subset{\Bbb R}^n$ is integration-by-parts
formulas, Green's formulas.
It is by no means obvious how one can establish such
formulas for the present nonlocal operators. Interesting generalizations have recently been obtained for 
translation-invariant operators
 by Ros-Oton and Serra,  partly with Valdinoci,
 in \cite{RS14b,RSV15}, and applied to nonlinear equations $Pu=f(u)$
 there as well as in \cite{RS14c,RS15a}; they have also been applied to
 nonlinear time-dependent
 Schr\"odinger equations by Boulenger, Himmelsbach and Lenzmann in \cite{BHL16}.

In the present paper we show an extension of the formulas to
$x$-dependent pseudodifferential operators, by completely different
methods.
The key ingredient is a factorization of $P$, $P\sim P^-P^+$ modulo
certain smoothing operators, where $P^-$ and $P^+$ preserve support in
$\complement\Omega $ resp.\  $\overline\Omega $. The operators $P^\pm$
are not standard $\psi $do's; their symbol structure is analyzed in
detail, using symbol classes originally introduced for Poisson and
trace operators in the Boutet de Monvel calculus \cite{B71}, as
developed in \cite{G90,G96}.

Our main results are: When $P$ is a classical  $x$-dependent $\psi
$do of order $2a$ on ${\Bbb R}^n$ with even symbol, elliptic avoiding a ray, and $\Omega $ is a
smooth bounded subset of $\R^n$, then the solutions $u$ of the Dirichlet
problem
$$
r^+Pu=f\text{ on }\Omega ,\quad \operatorname{supp}u\subset\comega,\tag1.3
$$
satisfy
$$
\align
\int_{\Omega  }( Pu\,\partial_j\bar
u'+\partial_ju\,\overline{  P^*u'})\,dx
&=\Gamma
(a+1)^2\int_{\partial\Omega }\nu _js_0\gamma _0(d^{-a}u)\,\gamma
_0(d^{-a}\bar u')\, d\sigma \\
&+\int_{\Omega } [P,\partial_j]u\,\bar u'\,dx,\quad j=1,\dots,n.\tag1.4
\endalign
$$
Here 
$\nu =(\nu _1,\dots,\nu _n)$ is the interior normal vector field to
$\partial\Omega $, and $s_0(x)$ is the principal symbol of $P$ at
$(x,\nu (x))$. As a corollary, we find
$$
\align
\int_{\Omega  }( Pu\,(x\cdot \nabla \bar u')+(x\cdot
\nabla u)\, \overline {P^*u'})\,dx
&=\Gamma
(a+1)^2\int_{\partial\Omega }(x\cdot\nu)s_0 \gamma _0(d^{-a}u)\,\gamma
_0(d^{-a}\bar u')\,
d\sigma\\
&-n\int_{\Omega } Pu\,\bar u'\,dx+\int _\Omega [P,x\cdot\nabla]u\,\bar
u'\,dx;
\tag1.5
\endalign
$$
where in detail,
$$
[P,x\cdot\nabla]=P_1    -P_2, \quad P_1=\Op(\xi \cdot\nabla_\xi
p(x,\xi )),\quad P_2=\Op(x
\cdot\nabla_x p(x,\xi )).\tag1.6
$$
 For $P$ one can for example take $A(x,D)^a$, where
$A(x,D)$ is a strongly
elliptic  second-order differential operator, in particular e.g.\ $(-\Delta +A_1(x,D))^a$, where $A_1$ is of order 1.

The formulas hold when
$u$ and $u'$ are solutions of problems (1.3) with  
$f\in \ol H^{1-a }(\Omega )$, or $f\in C^{1-a+\varepsilon }(\comega)$
for some $\varepsilon >0$. They extend to  $f\in \ol H^{\frac12-a +\varepsilon
}(\Omega )$, when the integrals are understood as Sobolev space dualities.

The formulas shown in \cite{RS14b,RSV15} that (1.4)--(1.5) extend, are concerned with real solutions
of (1.3) with $f\in C^{0,1}(\comega)$; here $\Omega $ is a bounded
$C^{1,1}$-domain, and $P$ is translation-invariant with even nonnegative
homogeneous kernel,
with possibly less smoothness than in our $C^\infty $-case.
In comparison, our method allows nonselfadjointness, nonpositivity, $x$-dependence, and
nonhomogeneity (lower-order terms).

To our knowledge, this is the first time that the integration-by-parts
question with $P$ and $\partial_j$ (or $x\cdot\nabla$) has  been solved for operators that
are not translation-invariant. Note that the $x$-dependence results in a
new interior term with $[P,\partial _j]$ in (1.4) (and with
$P_2$ in (1.5)--(1.6)).

The version of (1.5) shown in  \cite{RS14b,RSV15} is important in the study of existence questions for
nonlinear problems  where $f$ is replaced by $f(u)$ in (1.3), since it
leads to a Pohozaev-type identity for the possible
solutions. Here we find the following extended Pohozaev identity, for 
selfadjoint $x$-dependent nonhomogeneous operators $P$:
$$
\aligned
-2n\int_\Omega F(u)\,dx+n\int_\Omega
f( u)\,u\, dx&=\Gamma (1+a)^2\int_{\partial\Omega }(x\cdot \nu) \,s_0
\gamma _0(d^{-a}u)^2
\, d\sigma \\
&\quad +\int_{\Omega } [P,x\cdot\nabla]u\,u\,dx, 
\endaligned
\tag1.7
$$ 
for real solutions $u$; here
$F(t)=\int_0^tf(s)\,ds$. As a simple example, we apply the formula to 
$P=(-\Delta +m^2)^a$, $m>0$, showing a unique continuation
principle, and nonexistence of bounded nontrivial solutions to (1.3)
with $f$ taken as $f(u)=\operatorname{sign}u\,|u|^{r}$, $r\ge \frac{n+2a}{n-2a}$.
\medskip

{\it Plan of the paper:} The Appendix contains the notation for
function spaces, and collects some facts on pseudodifferential operators
that are known from the general theory and from preceding works such as
\cite{G15,G14}. Section 2 shows the factorization of symbols having
the $a$-transmission property, and describes the symbol spaces and
mapping properties of the generalized $\psi $do's that arise from the
construction. In Section 3 we establish the formula (1.4) in the case
where $\Omega $ is replaced by $\rnp$, for $j=n$. Finally
in Section 4, we treat the problem for arbitrary smooth domains
$\Omega $, showing the formulas (1.4)--(1.7) in general and drawing
some consequences.

 \head 2. Factorization of homogeneous symbols \endhead

\subhead{2.1 Some notation}\endsubhead

The function $\ang\xi $ stands for $(1+|\xi |^2)^\frac12 $, and we denote by $[\xi ]$ a positive
$C^\infty $-function equal to $|\xi |$ for $|\xi |\ge 1$ and $\ge
\tfrac12$ for all $\xi $. Multi-index notation is used for
differentiation (and polynomials):
$\partial=(\partial_1,\dots,\partial_n)$, and $\partial^\alpha
=\partial_1^{\alpha _1}\dots \partial_n^{\alpha _n}$ for $\alpha
\in{\Bbb N}_0^n$, with $|\alpha |=\alpha _1+\cdots+\alpha _n$, $\alpha
!=\alpha _1!\dots\alpha _n!$. $D=(D_1,\dots,D_n)$ with $D_j=-i\partial_j$.

Operators are considered acting on functions or distributions on
${\Bbb R}^n$, and on subsets  
 $\rnpm=\{x\in
{\Bbb R}^n\mid x_n\gtrless 0\}$ (where $(x_1,\dots, x_{n-1})=x'$), and
 bounded $C^\infty $-subsets $\Omega $ with  boundary $\partial\Omega $, and
their complements.

Restriction from $\R^n$ to $\rnpm$ (or from
${\Bbb R}^n$ to $\Omega $ resp.\ $\complement\comega$) is denoted $r^\pm$,
 extension by zero from $\rnpm$ to $\R^n$ (or from $\Omega $ resp.\
 $\complement\comega$ to ${\Bbb R}^n$) is denoted $e^\pm$. Restriction
 from $\crnp$ or $\comega$ to $\partial\rnp$ resp.\ $\partial\Omega $
 is denoted $\gamma _0$. 

The reader
 is encouraged to consult the Appendix for further notation, as it
 becomes relevant.

\subhead{2.2 The factorization question}\endsubhead 

Let there be given a function $p(\xi )\in C^\infty ({\Bbb R}^n\setminus \{0\})$, homogeneous of degree $2a$ with $0<a<1$, even and
elliptic, i.e.,
$$
p(-\xi )=p(\xi )\text{ and } p(\xi )\ne 0 \text{ for all }\xi \in{\Bbb
R}^n\setminus \{0\}.\tag2.1
$$

Consider the points in ${\Bbb R}^n$ as $\xi =(\xi ',\xi _n)$, where
$\xi '\in{\Bbb R}^{n-1}$, $\xi _n\in{\Bbb R}$. 
According to Vishik and Eskin, see Eskin \cite{E81} Ch.\ 6, $p$ can be
written as a product of two factors $p_+(\xi ',\xi _n)$ and $p_-(\xi ',\xi
_n)$ that extend analytically in $\xi _n$
to ${\Bbb C}_-$ resp.\ ${\Bbb C}_+$; here ${\Bbb C}_{\pm}=\{\xi
_n\in{\Bbb C}\mid \operatorname{Im}\xi _n\gtrless 0\}$.

Since the sign convention for the Fourier transform in \cite{E81} is
the opposite of the standard choice in Western literature, with
consequences for other $\pm$-conventions, it is hard to avoid
confusion when quoting the book directly. Therefore we shall show a 
detailed version of the factorization, where we moreover relate it to  the
symbol estimates and points of view that play a role in \cite{H65},
\cite{B71} and later works such as \cite{G96}, \cite{G09}, \cite{G15}.

By division by the number $p(0,1)$, we can assume that $p(0,1)=1$.

We define (for $\xi \ne 0$)
$$
q(\xi )=p(\xi )|\xi |^{-2a}, \quad \psi (\xi )=\operatorname{log}q(\xi ),\tag2.2
$$
they are both homogeneous of degree 0 and even. Actually, it suffices
for the following considerations that $q$ is ``even in the $\xi
_n$-direction'', more precisely, has the 0-transmission property
(\cite{B71,G96,G09,G15}) with
respect to the surfaces $\{x_n=c\}$:
$$
\partial_\xi ^\alpha q(0,-\xi _n)=(-1)^{|\alpha |}\partial_\xi
^\alpha q(0,\xi _n),\text{ all }\alpha \in{\Bbb N}_0^n,\tag2.3
$$
which clearly holds for even symbols of order 0. In order to have
the logarithm defined bijectively, we assume that the values of  $q$ 
avoid some ray $\{z=re^{i\theta} \mid r\ge 0\}$ in ${\Bbb
C}$. (\cite{E81} includes some more
general symbols.)

When
$\xi '\ne 0$, 
$$
\lim_{\xi _n\to \pm \infty }q(\xi ',\xi _n)=\lim_{\xi _n\to \pm \infty }q(\xi '/\xi _n,1)
=1,\quad \lim_{\xi _n\to
\pm \infty }\psi (\xi ',\xi _n)=0.
\tag2.4
$$
To factorize $q$ we shall decompose $\psi $ into a sum of two terms that extend
holomorphically into ${\Bbb C}_\pm$, respectively. This can be
formulated in terms of Cauchy integral formulas. 

Let us recall some facts about Cauchy integral decompositions. 
When $f(t )$ is $O(\ang t ^{-1})$ on ${\Bbb R}$ with a
continuous  derivative
$f'(t )$ that is $O(\ang t ^{-2})$ on ${\Bbb R}$, one can define
$$\aligned
f _+(t)&=\frac {i}{2\pi}\int _{{\Bbb R}} \frac{f(\sigma )}{\sigma
-t }\,\d\sigma \text{ for }\operatorname{Im}t<0,\\
f _-(t)&=\frac {-i}{2\pi}\int _{{\Bbb R}} \frac{f(\sigma )}{\sigma
-t }\,\d\sigma \text{ for }\operatorname{Im}t>0;
\endaligned\tag2.5
$$
they are holomorphic for $t\in{\Bbb C}_-$ resp.\ ${\Bbb C}_+$, and extend by
continuity to  $\overline{\Bbb C}_-$ resp.\ $\overline{\Bbb C}_+$ The
values on ${\Bbb R}$ (the limits for $\operatorname{Im}t\to 0$ from
${\Bbb C}_-$ resp.\ ${\Bbb C}_+$) satisfy
$$
f_+(t )+f_-(t )=f(t ).\tag2.6
$$
Moreover, for the functions on ${\Bbb R}$, the inverse Fourier 
transforms satisfy
$$
\Cal F^{-1}f_{+}=e^+r^+\Cal F^{-1}f,\quad \Cal F^{-1}f_{-}=e^-r^-\Cal F^{-1}f.\tag2.7
$$
(they are in $L_2({\Bbb R})$); here $r^\pm$ denotes restriction from
functions on ${\Bbb R}$ to functions on ${\Bbb R}_\pm$, and $e^\pm$ denotes
extension of functions on ${\Bbb R}_\pm$ to functions on ${\Bbb R}$ by zero on ${\Bbb R}_\mp$.
These facts are well-known; proofs can be found e.g.\ in \cite{E81}
Lemma 6.1, Th.\ 5.1. (\cite{H65} refers for the decomposition to
Beurling's contribution to the Helsingfors congress 1938.)

As in \cite{B71,G96,G09} we shall denote the mappings by $h^\pm\colon f\to
f_\pm$; note that they are complementing projections, satisfying $h^++h^-=I$.
(The mappings $h^\pm$ correspond to the mappings $\Pi ^\pm$ in \cite{E81}, except that the holomorphy regions are
exchanged because of a different convention for the Fourier transform.)
The mappings are applied to special spaces of $C^\infty $-functions in
the calculus of \cite{B71}; there are detailed accounts e.g.\ in
\cite{G96} Sect.\ 2.2 or \cite{G09}  Ch.\ 10, which serve our
purposes here (and will be taken up below in Section 2.3). The projection properties are summed up e.g.\ in \cite{G09} Th.\ 10.15. 

Recall some
ingredients: 
With $d\in{\Bbb Z}$, 
$\Cal H_{d}$ denotes the space of $C^\infty $-functions $f(t)$
on ${\Bbb R}$ such that $k(\tau )=\tau ^df(\tau ^{-1})$ coincides with
a $C^\infty $-function for $-1<\tau <1$ (this means that the
derivatives of $f$ match in a good way for $t\to\pm\infty $). Here one
can show that 
$$
\Cal H_{-1}=\Cal F(e^-\Cal S_-\oplus e^+\Cal S_+),\quad \Cal
H_{d}=\Cal H_{-1}\oplus {\Bbb C}_d[t] \text{ for }d\ge 0,\tag2.8
$$
where $\Cal S_{\pm}=r^\pm \Cal S({\Bbb R})=\Cal S(\crpm)$ (defined
from the Schwartz space $\Cal S({\Bbb R})$), and
${\Bbb C}_d[t]$ stands for the space of polynomials of degree $\le d$
in $t$.
Setting (with a slight asymmetry)
$$
\Cal H^+=\Cal F (e^+\Cal S_+), \quad \Cal H^-_{d}=\Cal F (e^-\Cal
S_-)\oplus {\Bbb C}_d[t],  \tag2.9
$$
one defines the mappings $h^\pm$ on $\Cal H_d$, consistently with their definition
given above for $d\le -1$, such that they are projections with ranges 
$$
h^+\Cal H_d=\Cal H^+,\quad h^-\Cal H_d=\Cal H^-_d ,\quad\text{for
}d\ge -1.\tag2.10
$$
The symbol $h_{-1}$ denotes the projection from $\Cal H_d$ to $\Cal
H_{-1}$ that removes the polynomial part.
The space $\Cal H^-_{-1}$ equals the space of conjugates of functions
in $\Cal H^+$ (\cite{G09} (10.55)). $\Cal H^+$ can also be denoted
$\Cal H^+_{-1}$ when convenient. Note that when $f\in\Cal H_{-1}$,
$\overline{h^-f}=h^+(\overline f)$.

In the case we shall work on, 
we are looking for a {\it factorization}, not a sum decomposition.

This was not treated in \cite{B71, G09}. It involves taking the
logarithm of $q$, decomposing $\log q$ into a sum by Cauchy integrals, and then deriving a
factorization of $q$ itself by exponentiating. The method is described
in \cite{E81} with a few estimates, but it has not been worked out what
happens in terms of $\Cal H^\pm$ spaces, so a new analysis is needed
for our purposes. Here we moreover find a special structure of the
factors, that in our application later will allow an
integration by parts formula. 

We first introduce some generalized symbol spaces and $\psi $do's.

\subhead{2.3 Symbol spaces for generalized $\psi $do's}\endsubhead

Homogeneous functions of $\xi $ are usually singular at $\xi =0$. We
use in general the convention that a symbol $p(x,\xi )$ is assumed to
be $C^\infty $ for all $\xi $, then in the homogeneous case,
homogeneity is assumed only for
$|\xi |\ge 1$, or $|\xi |\ge \delta $ for a suitable $\delta
>0$ (if needed, the associated fully homogeneous function is then
called the {\it strictly} homogeneous symbol). 

{\it Classical (also called polyhomogeneous) $\psi $do symbols of order $m$} are $C^\infty
$-functions having asymptotic series expansions
$p(x,\xi )\sim \sum_{j\in{\Bbb N}_0}p_j(x,\xi )$, where the $p_j$ are
homogeneous of degree $m-j$ in $\xi $ for $|\xi |\ge 1$ and 
$
\partial_x^\beta \partial_\xi ^\alpha (p(x,\xi )-\sum_{j<J}p_j(x,\xi
))\text{ is }O(\ang\xi ^{m-J-|\alpha |})
$ for all $\alpha ,\beta ,J$.
  The replacement of a strictly homogeneous function by a function
  that is smooth near $\xi =0$ is often achieved by multiplication by
  an {\it excision function $\eta (\xi )$} satisfying:
$$
\eta (\xi )=\eta (|\xi |)\in C^\infty ({\Bbb R}^n, [0,1])\text{ with } \eta (\xi )=0\text{
for }|\xi |\le \tfrac12, \; \eta (\xi )=1\text{
for }|\xi |\ge 1.\tag2.11
$$

It is a basic fact in the Boutet de Monvel calculus (cf.\ e.g.\ \cite{G09}
Th.\ 10.21) that when
$q(x,\xi )$ is a $\psi $do symbol of order $d\in{\Bbb Z}$ having the
0-transmission property with respect to the hyperplane  $\{x _n=c\}$, then the symbol $q(x',c,\xi ',\xi _n )$ is in $\Cal
H_d$ as a function of $\xi _n$, and 
$$
h^+q(x',c,\xi ',\xi _n )\in S^d_{1,0}({\Bbb R}^{n-1}, {\Bbb R}^{n-1}, \Cal H^+),\tag2.12
$$
where $h^+\colon f\mapsto f_+$ is the projection defined in (2.5)ff.\
(the space $S^d_{1,0}({\Bbb R}^{n-1}, {\Bbb R}^{n-1}, \Cal H^+)$ will be recalled in a moment). The function $h^+q$
is not quite a $\psi $do symbol in $\xi $ (although it is so in $\xi
'$ for each $\xi _{n}$), but we can still use the $\Op$-definition
(as in (A.1)),
and we call such symbols generalized $\psi $do symbols.

The symbol
spaces are explained e.g.\  in \cite{G09}, Section 10.3. With $m$
denoting a positive integer,
$
S^d_{1,0}({\Bbb R}^{m},{\Bbb R}^{n-1},\Cal H^+)$ consists of 
the following $C^\infty $-functions:
$$
\aligned
f(X,\xi ',\xi _n)&\in S^d_{1,0}({\Bbb R}^{m},{\Bbb R}^{n-1},\Cal
H^+),\text{ when }f(X,\xi ',\xi _n)\text{ is in }\Cal H^+\text{ w.r.t.\
}\xi _n\text{, and} 
\\
\| D_{X}^\beta D_{\xi '}^\alpha D_{\xi _n}^k&h_{-1}(\xi _n^{k'}f(X,\xi
',\xi _n))\|_{L_2({\Bbb R})}\le C_{\alpha ,\beta ,k,k'} \ang{\xi
'}^{d+\frac12-k+k'-|\alpha |}, 
\endaligned\tag2.13
$$
for all indices $\alpha \in{\Bbb N}_0^{n-1},\beta \in{\Bbb N}_0^m,k,k'\in{\Bbb N}_0$, with constants $
C_{\alpha ,\beta ,k,k'} $. $m$ is usually taken equal to $n$ or $n-1$. 
(The definition in \cite{G09} has $h^+$
instead of $h_{-1}$; the projections $h^+$ and $h_{-1}$ have the same
effect of removing
the polynomial terms arising from the multiplication of an $\Cal
H^+$-function by $\xi _n^{k'}$.)

The $L_2$-norms are useful when Fourier transforms are involved. In
fact, the system of seminorms (2.13) is {\it equivalent with} the following
system, applied to the inverse Fourier transforms  $\tilde f=\Cal
F^{-1}_{\xi _n\to x_n}f$ restricted to
$\{x_n>0\}$:
$$
\| D_{X}^\beta D_{\xi '}^\alpha x_n^kD_{x_n}^{k'}r^+\tilde f(X,x_n,\xi ')\|_{L_2(\rp)}\le C_{\alpha ,\beta ,k,k'} \ang{\xi '}^{d+\frac12-k+k'-|\alpha |}, \tag2.14
$$
the space of such functions $r^+\tilde f$ is denoted $S^d_{1,0}({\Bbb R}^{m},
{\Bbb R}^{n-1}, \Cal S_+)$. Here $\tilde f$ is in $e^+\Cal S_+$ as a
function of $x_n$. The effect of $h_{-1}$ is here replaced by that of
$r^+$, which removes possible linear combinations of $D^j_{x_n}\delta
_{x_n}$ (supported at $\{x_n=0\}$) that arise from differentiating
$\tilde f\in e^+\Cal S_+$.

It will be useful to observe that one can replace $L_2(\rp)$-norms
by $L_\infty (\rp)$-norms  or $L_1(\rp)$-norms (as remarked for $L_\infty
$-norms around (10.17) in \cite{G09}, and used sporadically in the
literature):

\proclaim{Lemma 2.1}
The family of estimates {\rm (2.14)} is equivalent with  the
family of estimates:
$$
\| D_{X}^\beta D_{\xi '}^\alpha x_n^kD_{x_n}^{k'}r^+\tilde f(X,x_n,\xi ')))\|_{L_\infty (\rp)}\le C_{\alpha ,\beta ,k,k'} \ang{\xi '}^{d+1-k+k'-|\alpha |}, \tag2.15
$$
as well as with the family of estimates
$$
\| D_{X}^\beta D_{\xi '}^\alpha x_n^kD_{x_n}^{k'}r^+\tilde f(X,x_n,\xi ')\|_{L_1(\rp)}\le C_{\alpha ,\beta ,k,k'} \ang{\xi '}^{d-k+k'-|\alpha |}. \tag2.16
$$
\endproclaim 

\demo{Proof}
We have the elementary
inequalities for functions $u(t)\in\Cal S_+$, $\sigma >0$:
$$
\aligned
\sup_{t\ge 0}|u(t)|^2&\le \sup_{t\ge 0} \int_t^\infty
|\partial_{s}(u(s)\bar u(s))|\,ds\le
2\|u\|_{L_2(\rp)}\|\partial_tu\|_{L_2(\rp)},\\
\sup_{t\ge 0}|u(t)|&\le \sup_{t\ge 0} \int_t^\infty |\partial_{s}u(s)|\,ds\le
\|\partial_tu\|_{L_1(\rp)},\\
\|u\|_{L_2}&\le\big\|\frac {1+\sigma t}{1+\sigma t}u\big\|_{L_2}\le
c\sigma ^{-\frac12}\|(1+\sigma t)u\|_{L_\infty },
\\
\|u\|_{L_1}&=\int_0^\infty \frac {1+\sigma t}{1+\sigma
t}|u(t)|\,dt\le 
c\sigma ^{-\frac12}\|(1+\sigma t)u\|_{L_2},
\endaligned\tag2.17
$$
where $\|(1+\sigma t)^{-1}\|_{L_2}=c\sigma ^{-\frac12}$.

Thus when $u$ satisfies
$$
\| t^kD_{t}^{k'}u(t)\|_{L_2(\rp)}\le C_{k,k'}\sigma
^{d+\frac12-k+k'},\text{ all }k,k'\in{\Bbb N}_0, 
$$
then we have from the first line:
$$
\| u(t)\|_{L_\infty (\rp)}\le  (2C_{0,0}\,\sigma ^{d+\frac12}\,C_{1,0}\,\sigma ^{d+\frac32})^{\frac12}=c'\sigma ^{d+1}, 
$$
with a similar treatment of derived functions $t^kD_t^{k'}u$. The variables $X,\xi '$ are easily
included, to see with $\sigma =\ang{\xi '}$ that the system of
estimates (2.14) implies (2.15). For the opposite direction, the basic
step is that when inequalities
$$
\|t^kD_t^{k'}u(t)\|_{L_\infty (\rp)}\le C_{k,k'}\sigma ^{d+1+k-k'}
$$
hold, then we have from the third line in (2.17) that
$$
\| u(t)\|_{L_2 (\rp)}\le c\sigma ^{-\frac12}(C_{0,0}\,\sigma ^{d+1}+\sigma\, C_{0,1}\, \sigma ^{d})=c''\sigma ^{d+\frac12}, 
$$
with a similar treatment of derived functions.

For $L_1$-norms, we moreover use the other lines in (2.17).
\qed
\enddemo

Instead of the above estimates that are global in $X$, we can work
with the constants $C_{\alpha ,\dots}$ replaced by continuous, hence locally bounded,
coefficients $C_{\alpha ,\dots}(X)$; they can be applied in localized
situations, and are more
general than the above.
Global estimates were considered in
\cite{G96,G09}, and are useful when one considers operators defined over
unbounded domains such as ${\Bbb R}^n$, $\rnp$ (more generally:
``admissible manifolds'', as defined in \cite{G96}).

We also need a notation for the spaces where the functions are in
$\Cal H^-_{-1}$ or in  $\Cal H_{-1}$ as functions of $\xi _n$: 
$$
\aligned
f(X,\xi ',\xi _n)&\in S^d_{1,0}({\Bbb R}^{m},{\Bbb R}^{n-1},\Cal
H^-_{-1}),\text{ when } f\in \Cal H^-_{-1}\text{ w.r.t.\ }\xi _n\text{
and}\\
\| D_{X}^\beta D_{\xi '}^\alpha D_{\xi _n}^k&h_{-1}(\xi _n^{k'}f(X,\xi
',\xi _n))\|_{L_2({\Bbb R})}\le C_{\alpha ,\beta ,k,k'} \ang{\xi
'}^{d+\frac12-k+k'-|\alpha |}, \\
f(X,\xi ',\xi _n)&\in S^d_{1,0}({\Bbb R}^{m},{\Bbb R}^{n-1},\Cal
H_{-1}),\text{ when } f\in \Cal H_{-1}\text{ w.r.t.\ }\xi _n\text{
and}\\
\| D_{X}^\beta D_{\xi '}^\alpha D_{\xi _n}^k&h_{-1}(\xi _n^{k'}f(X,\xi
',\xi _n))\|_{L_2({\Bbb R})}\le C_{\alpha ,\beta ,k,k'} \ang{\xi
'}^{d+\frac12-k+k'-|\alpha |}, 
\endaligned\tag2.18
$$
for all indices. Again, the estimates are equivalent with families of
estimates of the inverse Fourier transforms in $\xi _n$ as described
above for $\Cal H^+$. Note here that the inverse Fourier transform of
$\Cal H_{-1}=\Cal H^-_{-1}\oplus \Cal H^+$ is $e^-\Cal S_-\oplus
e^+\Cal S_+$, so that in fact, the second system of estimates is
equivalent with the system
$$
\| D_{X}^\beta D_{\xi '}^\alpha x_n^kD_{x_n}^{k'}\tilde f(X,x_n,\xi
')|_{\rmi\cup\rp}\|_{L_2(\rmi)\oplus L_2(\rp)}\le C_{\alpha ,\beta ,k,k'} \ang{\xi '}^{d+\frac12-k+k'-|\alpha |}. \tag2.19
$$

There are also versions of these spaces with local estimates in $X$ (i.e., with
the constants $C_{\alpha ,\dots}$ replaced
by continuous functions of $X$). 

The symbols in $S^d_{1,0}({\Bbb R}^{m},{\Bbb R}^{n-1},\Cal
H^+)$ were used in \cite{G96,G09} to define Poisson and trace
operators (maps between the boundary and the interior of $\rnp$). We shall here
use them to define $\psi $do's on ${\Bbb R}^n$. Since they do not satisfy
all the estimates usually required of $\psi $do symbols, we view them
as {\it generalized} $\psi $do symbols, and the operators resulting from
applying the Op-definition in (A.1) as {\it generalized} $\psi $do's.
To find their mapping properties, it is important to derive relevant sup-norm estimates
from (2.13) (and here it is a point to avoid having to involve the
projection $h_{-1}$). 

\proclaim{Lemma 2.2} Let $f\in S^d_{1,0}({\Bbb R}^{m},{\Bbb R}^{n-1},\Cal
H_{-1})$. 

$1^\circ$ Then also $\xi _n^kD_{\xi _n}^kf$ is in the space for all
$k\in{\Bbb N}_0$,
and
$$
|D_{X}^\beta D_{\xi '}^\alpha \xi _n^kD_{\xi _n}^kf(X,\xi
',\xi _n)|\le C_{\alpha ,\beta ,k} \ang{\xi '}^{d-|\alpha |},\tag2.20
$$
for all $\alpha ,\beta ,k$.

$2^\circ$ Moreover, $(\ang{\xi '}\pm i\xi _n)D_{\xi _n}f$ belongs to $S^d_{1,0}({\Bbb R}^{m},{\Bbb R}^{n-1},\Cal
H_{-1})$.
\endproclaim

\demo{Proof} When $\varphi (\xi _n)\in \Cal H_{-1}$, then so are $D_{\xi
_n}\varphi $ and $\xi _nD_{\xi _n}\varphi $; without going deeply into
the definition of $\Cal H_{-1}$ and $h_{-1}$ we can see this by observing that the inverse Fourier
transforms $-x_n\tilde \varphi (x_n)$ and $-D_{x_n}x_n\tilde\varphi
(x_n)$ are in $e^-\Cal S_-\oplus e^+\Cal S_+$ without distribution terms supported at
$x_n=0$.

For $1^\circ$ we iterate these considerations, seeing that also $\xi _n^kD_{\xi _n}^k\varphi $ and
$D_{\xi _n}\xi _n^kD_{\xi _n}^k\varphi $ are in $\Cal H^+$. The
estimates in (2.18) then show that when $f\in S^d_{1,0}({\Bbb R}^{m},{\Bbb R}^{n-1},\Cal
H_{-1})$, then 
$$
\aligned
\| D_{X}^\beta D_{\xi '}^\alpha \xi _n^kD_{\xi _n}^kf(X,\xi
',\xi _n))\|_{L_2({\Bbb R})}&\le C_{\alpha ,\beta ,k} \ang{\xi
'}^{d+\frac12-|\alpha |}, \\
\| D_{X}^\beta D_{\xi '}^\alpha D_{\xi _n}\xi _n^kD_{\xi _n}^kf(X,\xi
',\xi _n))\|_{L_2({\Bbb R})}&\le C'_{\alpha ,\beta ,k} \ang{\xi
'}^{d-\frac12-|\alpha |}.
\endaligned
$$
This implies (2.20) by the first line in (2.17), extended to
functions on ${\Bbb R}$. 

The other estimates
needed for the space $S^d_{1,0}({\Bbb R}^{m},{\Bbb R}^{n-1},\Cal
H_{-1})$ follow easily by carrying the inspection a little further. This
shows $1^\circ$, and $2^\circ$ follows by adding a similar inspection of $\ang{\xi
'}D_{\xi _n}f$.
 \qed\enddemo

We now investigate the mapping properties of the generalized $\psi
$do's defined from these symbols.
Here it will be convenient to refer to not only the
$H^s_p$-spaces recalled in the Appendix, but also spaces with a
different differentiability degree in the $x_n$-direction (used
e.g.\ in \cite{H63,G96,G09} for $p=2$):
$$
H^{s,t}_p({\Bbb R}^n)=\{u\in \Cal S'({\Bbb R}^n\mid \Cal
F^{-1}(\ang\xi ^s \ang{\xi '}^t\hat u)\in L_p({\Bbb R}^n\}=\Xi ^{-s}{\Xi '}^{-t}L_p({\Bbb R}^n),
$$
 where $\Xi ^t=\Op(\ang\xi ^t)$, ${\Xi '}^t=\Op(\ang{\xi '}^t)$.

To simplify the notation, we in the following abbreviate $S^d_{1,0}({\Bbb
R}^n,{\Bbb R}^{n-1},\Cal H^+)$ to $S^d(\Cal H^+)$, and similarly
with $\Cal H^-_{-1}$ and $\Cal H_{-1}$.

\proclaim{Proposition 2.3} Let $f(x,\xi ',\xi _n)\in S^d (\Cal H_{-1})$ for some $d\in{\Bbb R}$. Then
$F=\Op(f)$ is continuous
$$
F\colon H_p^{s,t}({\Bbb R}^n)\to H_p^{s,t-d}({\Bbb R}^n)\text{ for all }s,t\in{\Bbb R}.\tag2.21
$$
\endproclaim

\demo{Proof} Consider first the case $d=0$. By Lemma 2.2, we have that
$$
\ang{\xi '}^{|\alpha |}D_{\xi '}^\alpha \xi _n^kD_{\xi _n}^kD_x^\beta f\text{ is
bounded for all }\alpha \in {\Bbb N}_0^{n-1},\beta \in{\Bbb N}_0^n,
k\in {\Bbb N}_0.
$$
Then Lizorkin's criterion assures that $F\colon L_p({\Bbb R}^n)\to
L_p({\Bbb R}^n)$ is bounded; this shows (2.21) for $s=t=0$. The use of
Lizorkin's criterion is explained e.g.\ in \cite{GK93} around Th.\
1.6, with references.

Next, observe that 
$$
\|Fu\|_{H^1_p}\le c(\sum_{j=1}^n\|D_jFu\|_{L_p}+\|Fu\|_{L_p}),
$$
where $D_jFu=\Op(\xi _jf+D_{x_j}f)u=FD_ju+\Op(D_{x_j}f)u$. Here
$$
\|FD_ju\|_{L_p}\le c\|D_ju\|_{L_p}\le c'\|u\|_{H^{1}_p}
$$
by the preceding result, and since $D_{x_j}f$ is also in
$S^{d} (\Cal H_{-1})$,
$$
\|\Op(D_{x_j}f)u\|_{L_p}\le c\|u\|_{L_p}\le c'\|u\|_{H^{1}_p},
$$
implying altogether that
$$
F\colon H^{1}_p\to H^{1}_p
$$
is bounded. The argument extends easily to higher derivatives,
implying
boundedness of $$
F\colon H^{s}_p\to H^{s}_p \tag2.22
$$
for all $s\in {\Bbb N}_0$. By interpolation, the result extends to $s\in\crp$.

Since the Lizorkin criterion also holds when the operator is in
$y$-form (cf.\ \cite{GK93}), we likewise find that the adjoint
operator $F^*$ satisfies (2.22) for $s\ge 0$. Her we can replace $p$
by $p'$, and hence conclude by duality that (2.22) holds for $F$, all
$s\in {\Bbb R}$.

To extend the result to $H^{s,t}_p$-spaces, a first step is to observe that
$$
\|Fu\|_{H^{s,1}_p}\le
c(\sum_{j=1}^{n-1}\|D_jFu\|_{H^{s}_p}+\|Fu\|_{H^{s}_p})\le
c'(\sum_{j=1}^{n-1}\|D_ju \|_{H^s_p}+\|u\|_{H^s_p})\le c''\|u\|_{H^{s,1}_p}.
$$
Generalizing this to higher derivatives, we find (2.21) for
$s\in{\Bbb R}$, $ t\in{\Bbb N}_0$, $ d=0$, and interpolation and a similar
treatment of the adjoint leads to (2.21) for all $s,t\in{\Bbb R}$
when $d=0$.

For general $d\in{\Bbb R}$, we observe that $F{\Xi
'}^{-d}=\Op(f\ang{\xi '}^{-d})$, where $f\ang{\xi '}^{-d}\in
S^0 (\Cal H_{-1}) $, hence satisfies (2.21) with $d=0$. Then since obviously ${\Xi '}^d\colon
H^{s,t}_p({\Bbb R}^n)\simto H^{s,t-d}_p({\Bbb R}^n)$; (2.21) follows
for $F=F{\Xi '}^{-d}{\Xi '}^d$.
\qed
\enddemo

\proclaim{Theorem 2.4} Let $f(x,\xi )\in S^{d} (\Cal H_{-1})$
for some $d\in{\Bbb Z}$. Then $F=\Op(f)$ is continuous for all $s,t$:
$$
F\colon H^{s,t}_p({\Bbb R}^n)\to  H^{s-d,t}_p(\R^n)\text{ if
}d\ge -1.
\tag2.23
$$
The mapping property extends to $d=-k-1$, $k\in{\Bbb N}$, if $f(x,\xi
)([\xi ']+i\xi _n)^k\in  S^{-1} (\Cal H_{-1})$ (or  $f(x,\xi
)([\xi ']-i\xi _n)^k\in  S^{-1} (\Cal H_{-1})$).

\endproclaim

\demo{Proof} 
For $d\ge 0$, the result follows immediately from Proposition 2.3 
since $H^{s,t-d}_p({\Bbb
R}^n)\subset H^{s-d,t}_p({\Bbb R}^n)$.

Now let $d=-1$. 
Observe that
$$
F=F\,\Xi _+^{1}\, \Xi _+^{-1},\text{ where }F\,\Xi _+^{1}=\Op (f(x,\xi )([\xi ']+i\xi _{n})).
$$
This symbol is in $\Cal H_0$ as a function of $\xi _n$, and can be
decomposed as
$$
f(x,\xi )([\xi ']+i\xi _{n})=f(x,\xi )[\xi ']+h_{-1}(if\xi _{n})+(1-h_{-1})(if\xi _n).
$$
The first two terms are in $S^{0} (\Cal H_{-1})$, hence the
corresponding operators act as in (2.23) for $d=0$. The third term is
of the form $s(x,\xi ')$,  
constant in $\xi _n$ and with estimates $D_x^\beta D_{\xi }^\alpha
s=O(\ang{\xi '}^{-|\alpha |})$ (it is the
zero'th term in the expansion of $if\xi _n$ in powers $\xi _n^{-j}$,
$j\in{\Bbb N}_0$, cf.\ e.g.\ \cite{G09}, Def.\ 10.12).
It likewise defines a bounded operator in
$H^{s,t}_p({\Bbb R}^n)$.
Since $\Xi
_+^{-1}\colon H^{s,t}_p(\R^n)\simto H^{s+1,t}_p(\R^n)
$, we conclude (2.23) for $d=-1$. Note that we could just as well have
used compositions to the right with $\Xi _-^{\pm 1}=\Op(([\xi ']-i\xi
_n)^{\pm 1})$.

For the lower values of $d$ we apply the case $d=-1$ to the symbol
$f(x,\xi )([\xi ']+i\xi _n)^k$ (resp.\ $f(x,\xi )([\xi ']-i\xi _n)^k$). 
\qed
\enddemo

The most general symbols in $S^{-k-1} (\Cal H_{-1}$), $k\in{\Bbb
N}$, only have the
mapping property
$$
F\colon H^{s,t}_p({\Bbb R  }^n)\to  H^{s+1,t+k}_p(\R^n),
$$
(since they may only be $O(\xi _n^{-1})$ for $\xi _n\to\pm\infty $); this
is shown by combining (2.23) for $d=-1$ with Proposition 2.3. Fortunately, our applications in this paper will mainly be in the cases $d=0$ and
$d=-1$. Therefore we shall not burden the exposition with additional
terminology for symbol classes.

\subhead{2.4 The basic factorization theorem}\endsubhead 

With these
preparations, we shall etablish the factorization theorem for
homogeneous symbols.

\proclaim{Definition 2.5} When $P$ is a pseudodifferential operator
on $\R^n$ with a classical symbol $p(x,\xi )$ of order $m$, we say
that $P$ (and $p$) is {\bf elliptic avoiding a ray} $re^{i\theta  }$ when, for some $\theta
\in[0,2\pi ]$, the principal symbol  $p_0(x,\xi )$ takes values in ${\Bbb
C}\setminus\{z=re^{i\theta }\mid r\ge 0\} $ for all $\xi \in\R^n$ with
$|\xi |\ge 1$.
\endproclaim

The ray condition is usually assumed in resolvent constructions, and is
sometimes called "Agmon's condition", or "the condition for having a ray
of minimal growth". It is satisfied in particular when $P$ is {\it
strongly elliptic}, i.e., $\operatorname{Re}p_0(x,\xi )>0$ for $|\xi|
\ge 1$.

When the ray condition holds, $p_0(x,\xi )$ can be extended smoothly
into $|\xi |<1$ such that $p_0(x,\xi )\in{\Bbb
C}\setminus\{z=re^{i\theta }\mid r\ge 0\} $ for all $x,\xi $.
Then when we define the logarithm to be continuous on ${\Bbb
C}\setminus\{z=re^{i\theta }\mid r\ge 0\} $, $p_0$ can be retrieved
from its logarithm,    $p_0(x,\xi )=\exp\log p_0(x,\xi) $.

The condition will allow us to use the projections $h^+,h^-$ and the symbol
spaces introduced above in a simple way.

\proclaim{Theorem 2.6} Let $q(x,\xi ',\xi _n)$ be a pseudodifferential
symbol of order $0$, homogeneous in $\xi $ of degree $0$ for $|\xi
|\ge 1$ and having
the $0$-transmission property at all hyperplanes $\{x_n=c\}$,  
$$
\partial_x^\beta \partial_\xi ^\alpha q(x,0,-\xi _n)=(-1)^{|\alpha |}\partial_x^\beta \partial_\xi
^\alpha q(x,0,\xi _n),\text{ all }\alpha , \beta \in{\Bbb N}_0^n, |\xi |\ge 1.\tag2.24
$$
and elliptic avoiding a ray $re^{i\theta }$.
Assume that $q(x,\xi) $ avoids the ray also for $|\xi |<
1$.
 Denote $q(x,0,1)=s_0(x)$. 

Then $q$ has a factorization
$$
q(x,\xi )=s_0(x)q^-(x,\xi )q^+(x,\xi ),\tag2.25
$$
where 
$q^\pm(x,\xi ',\xi
_n)$ are invertible, and extend holomorphically into  respectively ${\Bbb C}_\mp$ as functions of $\xi _n$. Moreover,
$$
q^+(x,\xi )=1+f(x,\xi )\text{ with } f(x,\xi )\in S^0_{1,0}({\Bbb
R}^{n},{\Bbb R}^{n-1},\Cal H^+),\tag2.26
$$
homogeneous of degree $0$ in $\xi $ for $|\xi '|\ge 1$, and
$\overline{q^-}$ has these properties too.
\endproclaim

\demo{Proof} By division by $s_0(x)$ we can obtain that
$q(x,0,1)=1$, which will be assumed from now on. Fix $x=(x',c)$. We shall suppress the explicit mention of $x$, since
the estimates of derivatives in $x$ follow in a standard way when the
claims have been shown with respect to $\xi $ at each $x$.

Define $\psi (\xi )=\log
q(\xi )$ (to take real values when $q(\xi )$ is positive). The function $\psi $ is
likewise homogeneous of degree 0 for $|\xi |\ge 1$, hence is 
a $\psi$do symbol; moreover, it again has the 0-transmission property. Since $q(0,1)=1$, we have (2.4) for
all $\xi '$. In view of the 0-transmission property, 
$\psi $ is in $ \Cal H_{-1}$ as a function of $\xi _n$ for each $\xi '$, and
by \cite{G09} Th.\ 10.21,
$$
\psi _+=h^+\psi \in S^{0}_{1,0}({\Bbb R}^{n-1},{\Bbb R}^{n-1},\Cal H^+),\tag2.27 
$$
 and is homogeneous in $\xi $ of degree 0 when $|\xi '|\ge
 1$; it extends holomorphically into ${\Bbb C}_-$ as a function of
 $\xi _n$. Moreover, we can define
$$
\psi _-=h^-\psi 
;\text{ then } \psi =\psi _++\psi _-,\tag2.28
$$
and $\overline{\psi _-}$ is similar to $\psi _+$ (since it equals
$h^+\overline \psi $).

We now form
$$
q^+(\xi )=\exp (\psi _+(\xi ))=1+\psi _+(\xi )+\tfrac12 \psi _+(\xi )^2+\dots,\tag2.29
$$
and $q^-=\exp(\psi _-)$, then $q=q^-q^+$. We
have to show the estimates claimed in (2.26). Let
$$
f=q^+-1=\sum_{k=1}^\infty \frac1{k!}\psi _+^k.\tag2.30
$$
Instead of considering $f$ directly, consider the inverse Fourier transform
from $\xi _n$ to $z_n$
(restricted to $z_n\in{\Bbb R}_+$),
$$
\tilde f=\sum_{k=1}^\infty \frac1{k!}\tilde \psi _+^{*k}, \quad \tilde
\psi _+^{*k}=\tilde\psi _+*\dots*\tilde\psi _+, \text{ ($k$ factors)}.\tag2.31
$$

Observe that
$$
\|\tilde \psi _+*\tilde \psi _+\|_{L_\infty} \le \|\tilde\psi
_+\|_{L_1}\|\tilde \psi _+\|_{L_\infty },\;\dots \;,
\|\tilde \psi _+^{*k}\|_{L_\infty} \le \|\tilde\psi
_+\|^{k-1}_{L_1}\|\tilde \psi _+\|_{L_\infty },\tag2.32
$$
so that $\sum_{k=1}^\infty \frac1{k!}\|\tilde \psi _+\|^{k-1}_{L_1}\|\tilde
\psi _+\|_{L_\infty }$ is a majorising series for the series in
(2.31). Hence it converges in $L_\infty $-norm, and the limit satisfies
the estimate
$$
\|\tilde f\|_{L_\infty }\le \sum_{k=1}^\infty \frac1{k!}\|\tilde \psi _+\|_{L_1}^{k-1}\|\tilde
\psi _+\|_{L_\infty }=\|\tilde \psi _+\|_{L_1}^{-1}\|(\exp(\|\tilde
\psi _+\|_{L_1}-1)\|\tilde \psi _+\|_{L_\infty }.\tag 2.33
$$
Since $\tilde \psi _+$ satisfies the estimates in (2.14), (2.15) and (2.16) with
$d=0$ (in particular, $\|\tilde\psi _+\|_{L_1}$ is bounded in $\xi
'$), this shows that $\tilde f$ satisfies the first estimate in (2.15),
with $d=0$.

The estimates of $z_n$-derivatives follow in the same way, when we note that
$D_{z_n}$ just hits one factor; the one on which we impose the
$L_\infty $-norm. Multiplication by a power $z_n^l$ of $z_n$
corresponds for the Fourier transform $f$ to a derivative $D_{\xi
_n}^l$, for which there is a Leibniz formula. We carry this over
to the terms in $\tilde f$, seeing that it produces an expression where it hits at
most $l$ factors in the product $\tilde\psi _+^{*k}$. When $k\to\infty $, the estimates give $k-l $ factors $\|\tilde \psi
_+\|_{L_1}$ besides at most $l$ factors where specific
estimates of functions derived from $\tilde \psi _+$ are needed. This allows a majorising
sequence, leading to the desired estimate for $\tilde f$. Derivatives
with respect to $\xi '$ and $x$ are straightforward to include.

There is a similar analysis of $q^-$.
\qed

\enddemo

Note that the theorem shows that the considered symbols have
factorization index 0 (the homogeneity degree of $q^+$). Symbols of
order 0 without the ray condition can have other integer factorization
indices, see \cite{E81}, Ch.\ 6.

It is important in Theorem 2.6 that $q$ (after division by $s_0$) is not just factored into $q^+$ and
$q^-$ with the mentioned estimates, but that the first term in each of
the two factors
is $1$, besides
a term with a decrease in $\xi _n$. This will be very useful in the
applications.  

\subhead 2.5 Factorization of full symbols \endsubhead

For a general polyhomogeneous symbol that is elliptic and of type 0,
 the above can be extended to a factorization (in the sense of
operator composition or Leibniz products) respecting also lower-order
terms. Recall the composition rule for $\psi $do's:
$$
\Op(a)\Op(b)\sim \Op(a\# b), \text{ where } a\# b= \sum_{\alpha \in{\Bbb
N}_0^n}\partial_\xi ^\alpha a(x,\xi ) D_x^ab(x,\xi )/\alpha !.\tag2.34
$$
The last expression is often called the Leibniz product of $a$ and
$b$.

We now show that it is possible to refine the factorization from
Theorem 2.6, taking lower-order
terms into account.

\proclaim{Theorem 2.7} Let $Q$ be a classical  $\psi $do on
${\Bbb R}^n$ of order $0$, with
symbol  $q\sim\sum_{j\in {\Bbb
N}_0}q_j$ (where $q_j(x,\xi )$ is homogeneous of degree $-j$ in $\xi $ for $|\xi
|\ge 1$), elliptic avoiding a ray,
and having the
$0$-transmission property with respect to all hyperplanes $\{x_n=c\}$, i.e.,
$$
\partial_x^\beta \partial_\xi ^\alpha q_j(x,0,-\xi _n)=(-1)^{j-|\alpha |}\partial_x^\beta \partial_\xi ^\alpha q_j(x,0,\xi _n)\text{ for }j\in {\Bbb N}_0,\alpha ,\beta \in{\Bbb N}_0^n,\; |\xi |\ge 1.\tag2.35
$$
Denote
$q_0(x,0,1)=s_0(x)$. 

There exist two generalized pseudodifferential
symbols $q^{\pm}(x,\xi )\sim \sum_{j\in {\Bbb
N}_0}q^\pm_j(x,\xi ) $, with  $q^\pm_j(x,\xi )$ homogeneous of degree
$-j$ in $\xi $ for $| \xi '|\ge 1$, 
such that 
$$
q_0^+(x,\xi )=1+f(x,\xi )\text{ with } f(x,\xi )\in S^0(\Cal H^+);\quad q^+_j\in S^{-j}(\Cal
H^+),\text{ for }j>0,\tag2.36
$$
 $\overline{q^-}(x,\xi )$ has a similar form, and
$$
q\sim s_0\,q^-\# q^+,\tag2.37
$$
in the sense that for all $K$, the difference between $s_0^{-1}q$ and the expression formed
of the terms in $q^+$ and $q^-$ down to order $-K$, composed by the
Leibniz formula applied for $|\alpha |\le K$, is in $ S^{-K-1}(\Cal
H_{-1})$.

From the symbols $q^\pm$ we can define generalized
$\psi $do's $Q^\pm$, 
respectively;
 then $Q-s_0Q^-Q^+$ has symbol in  $S^{-\infty } (\Cal
H_{-1})=\bigcap _d S^{d } (\Cal
H_{-1})$.

The operator  $Q-s_0Q^-Q^+$
is smoothing in the sense that it maps  $H^{s,t}(\R^n)$ to
\linebreak
$H^{s+1,\infty }(\R^n)=\bigcap_t H^{s+1,t}(\R^n)$, for all $s$. 
\endproclaim

\demo{Proof} By multiplication by $s_0^{-1}$ we can assume that
$q_0(x,0,1)=1$. The principal parts $q^\pm_0$ of $q^\pm$ are defined by
application of Theorem 2.6 to $q_0$. Now we have to construct the
lower-order symbols. This goes inductively as follows:

Collecting the terms of order $-1$ in (2.37) (cf.\ (2.34)), we find that $q_1^\pm$
should satisfy:
$$
q_1=q_0^-q_1^++q_1^-q_0^++\sum_{k\le n}\partial_{\xi
_k}q_0^-D_{x_k}q^+_0.
$$


Dividing by $q_0$ and using that $q_0=q_0^-q_0^+$, we can rewrite this as
$$
\frac{q_1^+}{q_0^+}+\frac{q_1^-}{q_0^-}=\frac
{q_1}{q_0}-\frac1{q_0}\sum_{k\le n}\partial_{\xi
_k}q_0^-D_{x_k}q^+_0,\tag2.38
$$
where the right-hand side is already known.
By Theorem 2.6, the function $q_0^+$ is 1 plus a function in $\Cal H^+$ at each $(x,\xi
')$, and since it is nonvanishing, the inverse is likewise of the form
in (2.26). The same holds for $\overline{q^-_0}$. 
Moreover, $q_1$ being of order $-1$ and having the $0$-transmission
property implies that it is in $\Cal H_{-1}$ as a function of $\xi _n$.
Thus the right-hand
side  of (2.38) is in $\Cal H_{-1}$, and the
left-hand side expresses a decomposition in its $\Cal H^+$-part and
$\Cal H^-_{-1}$-part, for each $(x,\xi ')$. The decomposition is
unique, and one checks that the two terms satisfy the appropriate
estimates.

This shows the first step, and in the general step, one similarly
determines the two terms  ${q_k^+}/{q_0^+}$ and ${q_k^-}/{q_0^-}$ as
the components in $\Cal H^+$ and
$\Cal H^-_{-1}$ of an expression formed of the preceding symbol terms
of the relevant homogeneity:
$$
\frac{q_k^+}{q_0^+}+\frac{q_k^-}{q_0^-}=\frac
{q_k}{q_0}-\frac1{q_0}\sum_{j+|\beta |=k, j<k}\frac1{\beta !}\partial_{\xi
}^\beta q_j^-D_{x}^\beta q^+_j.\tag2.39
$$

There is a standard way to associate an exact symbol $q^\pm(x,\xi )$ with the
series \linebreak $\sum_{j\in {\Bbb
N}_0}q^\pm_j(x,\xi )$, namely, a convergent sum $q^\pm(x,\xi )=\sum_{j\in {\Bbb
N}_0}\eta (\xi /t_j)q^\pm_j(x,\xi )$, where $t_j\to \infty $
sufficiently rapidly (for $\eta (\xi )$, see (2.11)). Any other choice of a symbol with the given
asymptotic expansion differs from this by a symbol in $ S^{-\infty }_{1,0}({\Bbb R}^{n},{\Bbb R}^{n-1},\Cal
H^\pm_{-1})$. 

Then one finds by use of the Leibniz formula and
regrouping of homogeneous terms of the same order, that $Q-s_0Q^-Q^+$ is
a generalized $\psi $do with symbol in $S^{-\infty } (\Cal
H_{-1})$. The last statement follows from  Theorem
2.4 ff.
\qed 
\enddemo

When $q$ is {\it
even} in $\xi $, that is,
$$
q_j(x,-\xi )=(-1)^jq_j(x,\xi )\text{ for }|\xi |\ge 1, \text{ all }x,\tag2.40
$$
the property (2.35) holds in any coordinate system.

We furthermore observe the following property.

\proclaim{Theorem 2.8} For the symbol $q^+$ constructed in Theorem
{\rm 2.7}, there is a
parametrix symbol
$\tilde q^+$ with similar symbol properties, such that 
$$
q^+\#\tilde q^+\sim 1\sim \tilde q^+\# q^+,\tag2.41
$$
in the space consisting of symbols in $ S^0_{1,0}({\Bbb R}^{n},{\Bbb R}^{n-1},\Cal
H_{-1})$ plus functions of $x$ (constant in $\xi $). 

There is a
similar result for $q^-$.
\endproclaim
\demo{Proof} We apply the standard parametrix construction: With
$\tilde q_0^+=1/q_0^+$, we have that
$$
q^+\#\tilde q_0^+=1+\sum_{\beta \ne 0}\tfrac1{\beta !}\partial_{\xi
}^\beta q_0^+D_{x}^\beta \tilde q^+_0+(q^+-q_0^+)\# \tilde q_0^+\sim 1+r,\tag2.42
$$
where $r\in S^{-1}_{1,0}({\Bbb R}^{n},{\Bbb R}^{n-1},\Cal H^+)$ is
defined from a regrouping of the terms according to homogeneity.
Then, defining
$$
\tilde r\sim \sum_{k=1}^\infty(-1)^kr^{\# k}, 
$$
where $r^{\# k}\sim r\# r\#\dots \# r$ with $k$ factors, we find that
$ \tilde q^+_0\# (1+\tilde r)$ is a right parametrix symbol for
$q^+$. Similarly, there is a left parametrix symbol, and they are
seen to be equivalent. Thus we can take $\tilde q^+= \tilde q^+_0\# (1+\tilde
r)$, and it has the asserted properties.\qed 

\enddemo 


\example{Remark 2.9}
The constructions in Theorems 2.6 and 2.7 have been developed from
\cite{H65}, Theorem 2.6.3 ff.\ in combination with our use of function
spaces based on $\Cal H^\pm$ as in \cite{G96}, \cite{G09}. The
purpose in \cite{H65} was to construct an operator that solves, in the parametrix
sense, certain boundary problems for operators such as e.g.\
$P=\Xi _-^aQ\Xi _+^a$ with nonzero boundary data and 0 data in the
interior, generalizing (2.45) below.  
For $Q$ itself, with
$\widetilde Q^+=\Op(\tilde q^+)$, it can be shown that the operator $K_{\widetilde Q^+}\colon \varphi
(x')\to \widetilde Q^+(\varphi (x')\otimes \delta ({x_n}))$, which is
a Poisson operator in the Boutet de Monvel calculus,  satisfies:
$$
r^+QK_{\widetilde Q^+}\colon \E'({\Bbb R}^{n-1})\to C^\infty (\crnp),\quad 
\gamma _{-1,0}K_{\widetilde Q^+}-I\colon \E'({\Bbb R}^{n-1})\to C^\infty ({\Bbb R}^{n-1}) \tag2.43
$$
(i.e.\ $r^+Q K_{\widetilde Q^+}$ and $\gamma _{-1,0}K_{\widetilde Q^+}-I$ are smoothing operators); hence $K_{\widetilde Q^+}$
defines a parametrix solution operator to the problem
$$
r^+Qw=0,\quad \gamma _{-1,0}w=\varphi .\tag 2.44
$$
Here $\gamma _{-1,0}$ is a generalization of $\gamma _{\mu ,0}$ to low
values of $\mu $, defined but not studied in detail in \cite{G15}. Then for $P$ one
finds, 
 setting $w=\Xi _+^au$, that $\Xi _+^{-a}K_{\widetilde
 Q^+}$
defines a parametrix solution operator to the problem
$$
r^+Pu=0,\quad \gamma _{a-1,0}u=\varphi ;\tag2.45
$$
here  $\Xi _+^{-a}K_{\widetilde
 Q^+}$ can be regarded as a (generalized) Poisson operator of noninteger order.
The problem was also discussed in \cite{G15}, Th.\ 6.5; the present construction gives
a more direct information. We shall possibly take up the details in another
publication.
\endexample

\example{Remark 2.10} The factorization idea $P=P^-P^+$ with factors
having opposite support preservation properties could also be used in
the proof of \cite{RSV15} for one-dimensional operators on rays,
instead of their factorization of selfadjoint positive operators as
$P=P^\frac12P^\frac12$.
\endexample



\head 3. Integration by parts for 
operators on
the half-space\endhead

The reader is encouraged to consult the Appendix for notation.

Let $P$ be a classical $\psi $do on ${\Bbb R}^n$ of order $2a$ ($0<a<1$), having the
$a$-transmission property at the boundary of $\rnp$. Recall from
\cite{G15} Th.\ 4.2 that $r^+P$ maps the space $H^{a(s)}(\crnp)$ (cf.\
(A.9)) continuously into $\ol H^{s-2a}(\rnp)$ when $s>a-\frac12$. 

We wish to reduce the expression
$$
\int_{\rnp}Pu\,\partial_n\bar u'\,
dx+\int_{\rnp}\partial_nu\,\overline {P^* u'}\, dx\tag3.1
$$
for functions $u,u'\in H^{a(s)}(\crnp)$ for suitable $s$, to an integral over
the boundary of suitable boundary values,  supplied in the
$x_n$-dependent case
with an extra integral over $\rnp$. The fact that we integrate over
$\rnp$ implies a restriction $r^+$ on the integrands, that we
therefore need not mention explicitly in the formula.

The central argument will first be presented in a simple constant-coefficient case.

\proclaim{Theorem 3.1} $1^\circ$ Let $u,u'\in \E_a(\crnp)$ with compact
support in $\crnp$. Let $w=r^+\Xi _+^au$, $w'=r^+\Xi _+^au'$. Then
$$
\int _{\rnp}\Xi _-^ae^+w\,\partial_n\bar u'\, dx=
(\gamma _0w,\gamma _0w')_{L_2(\R^{n-1})}+(w,\partial_nw')_{L_2(\rnp)}
.\tag 3.2
$$

$2^\circ$ The formula extends to $u,u'\in H^{a(s)}(\crnp)$ for $s>a+\frac12$,
 with dualities:
$$
\multline
\ang{r^+\Xi _-^a e^+w,\partial_n 
u'}_{ \ol H^{\frac12-a+\varepsilon }(\rnp) ,\dot H^{a-\frac12-\varepsilon
}(\crnp)}\\
=(\gamma _0w,\gamma _0w')_{L_2(\R^{n-1})}+\ang{w,\partial_n 
w'}_{ \ol H^{\frac12-\varepsilon }(\rnp) ,\dot H^{-\frac12+\varepsilon
}(\crnp)},
\endmultline\tag 3.3
$$
for any $0<\varepsilon \le s-a-\frac12$ with $\varepsilon <1$.

$3^\circ$ Here, when $s\ge a+1$, then {\rm (3.3)} can be written in the form {\rm
(3.2)}, all ingredients being locally integrable functions.
\endproclaim

\demo{Proof} $1^\circ$. First let $u,u'\in \E_a(\crnp)$ with compact
support. Since $u\in H^{a(s)}(\crnp)$ for any large $s$,
$w=r^+\Xi _+^au\in C^\infty (\crnp)\cap \ol H^s(\rnp)$ for any 
$s$, with $u=\Xi _+^{-a}e^+w$ (cf.\ \cite{G15}, Propositions 1.7 and
4.1). Moreover, $r^+\Xi _-^ae^+w\in  C^\infty (\crnp)\cap \ol H^s(\rnp)$ for any 
$s$. There is similar information for
$u',w'$. 

Since $u\in \E_a(\crnp)$ with compact support, $\partial_nu\in
\E_{a-1}(\crnp)$ with compact support. Here $x_n^{a-1}$ is integrable
over compact sets. Altogether,
$r^+\Xi
_-^ae^+w\, \partial_n\bar u$ is on $\crnp$ the product of $x_n^{a-1}$ with a
compactly supported smooth function, so the integral is well-defined.

We can also observe that by the identification of $e^+\ol H^t(\rnp)$ and $\dot
H^t(\crnp)$ for $|t|<\frac12$, $e^+w'\in \dot H^{\frac12-\varepsilon
}(\rnp)$ for any $\varepsilon \in \,]0,1[\,$, 
 so
$$
\partial_n u'=\partial_n\Xi _+^{-a}e^+w'\in \partial_n\dot
H^{a+\frac12-\varepsilon }(\crnp)\subset\dot
H^{a-\frac12-\varepsilon }(\crnp).\tag 3.4
$$
Then since  $r^+\Xi _-^ae^+\Xi _+^au\in  \ol
H^{\frac12-a+\varepsilon }(\rnp)$, the integral may be written as the duality
$$
I=\ang{r^+\Xi _-^a e^+w,\partial_n 
u'}_{ \ol H^{\frac12-a+\varepsilon }(\rnp) ,\dot H^{a-\frac12-\varepsilon
}(\crnp)}.
$$
 
Now note that by (A.7), $r^+\Xi _-^ae^+\colon \ol
H^{\frac12+\varepsilon }(\rnp)\to \ol
H^{\frac12-a+\varepsilon }(\rnp)$
has the adjoint \linebreak$\Xi _+^a\colon \dot H^{a-\frac12-\varepsilon
}(\crnp)\to \dot H^{-\frac12-\varepsilon
}(\crnp)$. We can then continue the calculation of $I$ as follows:
$$
I=\ang{w,\Xi _+^a\partial_n 
u'}_{ \ol H^{\frac12+\varepsilon } ,\dot H^{-\frac12-\varepsilon
}}=\ang{w,\partial_n\Xi _+^a 
u'}_{ \ol H^{\frac12+\varepsilon } ,\dot H^{-\frac12-\varepsilon
}}=\ang{w,\partial_ne^+w'}_{ \ol H^{\frac12+\varepsilon } ,\dot H^{-\frac12-\varepsilon}}. 
$$
Here $w'$ itself is a nice function on $\crnp$, but the extension
$e^+w'$ to $\R^n$
has the jump $\gamma _0w'$ at $x_n=0$,
and there holds the formula
$$
\partial_ne^+w'=\gamma _0 w'\otimes \delta (x_n)+e^+\partial_nw'.\tag3.5
$$
where $\otimes$ indicates a product of distributions with respect to
different variables ($x'$ resp.\ $x_n$). It is a distribution version
of Green's formula (cf.\ e.g.\ \cite{G96}
(2.2.38)--(2.2.39)). Recall moreover from distribution
theory (cf.\ e.g.\ \cite{G09} p.\ 307) that the ``two-sided'' trace operator $\widetilde\gamma _0\colon
v(x)\mapsto \widetilde \gamma _0 v=v(x',0)$ has the mapping
$\widetilde\gamma _0^*\colon \varphi (x')\mapsto \varphi (x')\otimes
\delta (x_n)$ as adjoint, with continuity properties
$$
\widetilde \gamma _0\colon H^{\frac12+\varepsilon }({\Bbb R}^n)\to H^{\varepsilon }({\Bbb
R}^{n-1}),\quad
\widetilde \gamma _0^*\colon H^{-\varepsilon }({\Bbb R}^{n-1})\to H^{-\frac12-\varepsilon }({\Bbb
R}^{n}), \text{ for }\varepsilon >0. \tag3.6
$$
Here $\widetilde\gamma _0^*\varphi $ is supported in $\{x_n=0\}$,
 hence lies in $\dot H^{-\frac12-\varepsilon }(\crnp)$. We can then write
$$
\partial_ne^+w'=\widetilde \gamma _0^*(\gamma _0 w')+e^+\partial_nw'.\tag3.7
$$

Since $w\in \ol H^{\frac12+\varepsilon }(\rnp)$, it has an extension
$W\in H^{\frac12+\varepsilon }(\R^n)$ with
$w=r^+W$, and  $\gamma _0 w= \widetilde
\gamma _0W$.
Then
$$
\ang{w,\widetilde \gamma _0^*(\gamma _0  w')}_{ \ol
H^{\frac12+\varepsilon }(\rnp), \dot H^{-\frac12-\varepsilon
}(\crnp)}=\ang{W,\widetilde \gamma _0^*(\gamma _0  w')}_{  H^{\frac12+\varepsilon }({\Bbb R}^n),  H^{-\frac12-\varepsilon }({\Bbb R}^n)};
$$
this is verified e.g.\ by approximating $\widetilde \gamma _0^*(\gamma
_0w')$ in $\dot H^{-\frac12-\varepsilon }$-norm by a sequence of
functions in $C_0^\infty (\rnp)$.  Here we can use (3.6) to write
$$
\multline
\ang{W,\widetilde \gamma _0^*(\gamma
_0  w')}_{H^{\frac12+\varepsilon }({\Bbb R}^n), H^{-\frac12-\varepsilon
}({\Bbb R}^n) }
=\ang{\widetilde \gamma _0W,\gamma
_0  w'}_{H^{\varepsilon }({\Bbb R}^{n-1}), H^{-\varepsilon }({\Bbb
R}^{n-1}) }\\
=\ang{ \gamma _0 w,\gamma
_0  w'}_{H^{\varepsilon }({\Bbb R}^{n-1}), H^{-\varepsilon }({\Bbb
R}^{n-1}) }=(
\gamma _0w,\gamma _0 w')_{L_2({\Bbb R}^{n-1})}.
\endmultline
$$
In the last step we used  that since both $\gamma _0w$ and
$\gamma _0w'$ are in $H^\varepsilon ({\Bbb R}^{n-1})\subset L_2({\Bbb
R}^{n-1}) $, the duality over the boundary is
in fact an $L_2({\Bbb R}^{n-1})$-scalar product. 

 Then finally
$$
\aligned
I&=\ang{w,\partial_ne^+w'}_{ \ol H^{\frac12+\varepsilon } ,\dot
H^{-\frac12-\varepsilon}}=\ang{w,\widetilde \gamma _0^*(\gamma _0
w')+e^+\partial_nw'}_{ \ol H^{\frac12+\varepsilon } ,\dot
H^{-\frac12-\varepsilon}}\\
&=(\gamma _0w,\gamma _0 w')_{L_2(\R^{n-1})}+\ang{w,e^+\partial_nw'}_{ \ol
H^{\frac12+\varepsilon } ,\dot H^{-\frac12-\varepsilon}}\\
&=(\gamma _0w,\gamma _0 w')_{L_2(\R^{n-1})}+(w,e^+\partial_nw')_{L_2(\rnp)},
\endaligned 
$$
where we used that $w'\in \bigcap_s\ol H^s (\rnp)$. This
shows (3.2).

$2^\circ$.
If $u,u'\in H^{a(s )}(\crnp)$ with $s>a+\frac12$, they are in $\in
 H^{a(a+\frac12+\varepsilon )}(\crnp)$ for an $\varepsilon \in
 \,]0,1[\,$, $\varepsilon \le s-a-\frac12$, and
 then $w,w'\in \ol
H^{\frac12+\varepsilon }(\rnp)$ by definition. 
Moreover, by (A.11), 
$$
\aligned
 u&\in e^+x_n^a \ol
H^{\frac12+\varepsilon }(\rnp )+\dot H^{a+\frac12+\varepsilon
 }(\crnp),\text{ hence}\\
 \partial_nu&\in e^+x_n^{a-1} \ol
H^{\frac12+\varepsilon }(\rnp )+e^+x_n^{a} \ol
H^{-\frac12+\varepsilon }(\rnp )
+\dot H^{a-\frac12+\varepsilon
 }(\crnp).
\endaligned\tag3.8
$$
(Since $\gamma _0u=0$, there is no distribution term supported at
$\{x_n=0\}$.)
On the other hand, since $e^+w=\Xi _+^au\in e^+ \ol
H^{\frac12+\varepsilon }(\rnp)\subset \dot H^{\frac12-\varepsilon
'}(\crnp)$,
$u=\Xi _+^{-a}e^+w$ satisfies
$$
\aligned
u&\in \dot H^{a+\frac12-\varepsilon '}(\crnp),\text{
any }\varepsilon '>0,\text{ hence}\\
\partial_nu&\in \dot H^{a-\frac12-\varepsilon '}(\crnp).
\endaligned\tag3.9
$$
There is similar information for $u'$.

Here we can approximate
$u,u'$ in the norm of $H^{a(a+\frac12+\varepsilon )}(\crnp)$ by compactly supported elements $u_k,u'_k$ of $\E_a(\crnp)$
(cf.\ \cite{G15} Prop.\ 4.1). Then $w_k=r^+\Xi _+^au_k$ and
$w'_k=r^+\Xi _+^au'_k$
converge in $\ol H^{\frac12+\varepsilon }(\rnp)$, and in particular,
$\partial_nu'_k$ converges in $\dot H^{a-\frac12-\varepsilon }(\crnp)$
and $\partial_nw'_k$ converges in $\ol H^{-\frac12+\varepsilon }(\rnp)=\dot H^{-\frac12+\varepsilon
}(\crnp)$. This implies (3.3) by passage to the limit, proving $2^\circ$.

$3^\circ$. If $s\ge a+1$, then $w,w'\in \ol H^1(\rnp)$, so $\partial_nw'\in
L_2(\rnp)$, and $r^+\Xi _-^ae^+w
\in \ol H^{1-a}(\rnp)\subset
L_2(\rnp)$. Moreover, by (A.11),$$
\aligned u&\in e^+x_n^a \ol
H^{1}(\rnp )+\dot H^{a+1}(\crnp),\text{ hence since }\gamma _0u=0,\\
\partial_nu&\in
e^+x_n^{a-1}\ol H^1(\rnp )+e^+x_n^{a}L_2(\rnp )+\dot
H^{a}(\crnp);
\endaligned\tag3.10$$
so $\partial_nu,\partial_nu'$ are functions.
\qed
\enddemo

\example{Remark 3.2} The distributional formulation (3.5) of Green's formula
 has been an important ingredient in systematic studies of boundary
 value problems for many years, for example in the construction of the
 Calder\'o{}n projector by Seeley \cite{S66}, H\"ormander \cite{H66}, see
 also \cite{G09}, Ch.\ 11. In the case $a=1$, Theorem 3.1 is quite
 elementary and can be shown by reference to the usual formulation of  Green's formula
$$
\int_{\rnp} \partial_n v\, \bar v'\, dx+\int_{\rnp }v\,\partial_n \bar v'\,
dx
= - \int_{R^{n-1}}\gamma _0v\gamma _0\bar v'\,dx'.\tag3.11
$$
Let $a=1$. Note that $\Xi ^1_\pm=\Xi '\pm \partial_n$, where $\Xi
'=\Op(\ang{\xi '})$, it acts in the $x'$-variable only, and is selfadjoint. Let $s=2$ for definiteness; here
 $H^{1(2)}(\crnp)=e^+\ol H^2(\rnp)\cap \dot H^1(\crnp)$ (one may consult Example 1.6 in \cite{G15}). 

Consider $u,u'\in e^+\ol H^2(\rnp)\cap \dot H^1(\crnp) $, and let 
$w=r^+(\Xi '+\partial_n)u$, $w'=r^+(\Xi '+\partial_n)u'$; they lie in
$\ol H^1(\rnp)$. Denote moreover
$v'=r^+\partial_nu'=w'-r^+\Xi 'u'$. Then
$$
\aligned
I&\equiv ( r^+(\Xi '-\partial_n)w,
r^+\partial_n u')_{\rnp}=(\Xi 'w,
v')_{\rnp}-(r^+\partial_nw, v')_{\rnp}\\
&= (\Xi 'w,
v')_{\rnp}+ (w,r^+\partial_n v')_{\rnp}+(\gamma
_0w,\gamma _0 v')_{\R^{n-1}},
\endaligned$$
using (3.11). Here we note that since $\gamma _0u'=0$, $\gamma
   _0v'=\gamma _0w'$. Now
$$
\aligned
I&= (w,\Xi 'v'
)_{\rnp}+ (w,r^+\partial_n v')_{\rnp}+(\gamma
_0w,\gamma _0 w')_{\R^{n-1}}\\
&= (w,r^+(\Xi '+ \partial_n)
v')_{\rnp}+(\gamma
_0w,\gamma _0 w')_{\R^{n-1}}=( w,r^+(\Xi '+\partial_n)\partial_nu
')_{\rnp}+(\gamma
_0w,\gamma _0 w')_{\R^{n-1}}\\
&=( w,r^+{\partial_nw'})_{\rnp}+(\gamma
_0w,\gamma _0 w')_{\R^{n-1}},
\endaligned$$
showing (3.2) in this case.
\endexample

An immediate consequence of Theorem 3.1 is the following integration-by-parts result
for fractional  Helmholtz  operators:

\proclaim{Theorem 3.3} Let $u$ and $u'$ be as in Theorem {\rm 3.1}
 $1^\circ$ or $3^\circ$. Then one has for $m>0$:
$$
\multline
\int _{\rnp}(-\Delta +m^2)^au\,\partial_n\bar u'\, dx+\int
_{\rnp}\partial_nu\,
{(-\Delta +m^2)^a\bar u'}\, dx\\
=\Gamma (a+1)^2\int_{{\Bbb
R}^{n-1}}\gamma _0(x_n^{-a}u)\,\gamma _0(x_n^{-a}\bar u')\, dx'.
\endmultline
\tag 3.12
$$
If $u$ and $u'$ are as in Theorem {\rm 3.1} $2^\circ$, the formula holds
with dualities, for small $\varepsilon >0$,
$$
\multline
\ang{r^+(-\Delta +m^2)^au,\partial_n 
u'}_{ \ol H^{\frac12-a+\varepsilon } ,\dot H^{a-\frac12-\varepsilon
}}
+\ang{\partial_nu,r^+(-\Delta +m^2)^au'}_{\dot H^{a-\frac12-\varepsilon
}, \ol H^{\frac12-a+\varepsilon } }\\
=\Gamma (a+1)^2\int_{{\Bbb
R}^{n-1}}\gamma _0(x_n^{-a}u)\,\gamma _0(x_n^{-a}\bar u')\, dx'.
\endmultline\tag3.13
$$

\endproclaim

\demo{Proof} We have that 
$$
(-\Delta +m^2)^a=\Op((|\xi |^2+m^2)^a)=\Xi ^a_{m,-} \Xi ^a_{m,+},
\quad \Xi ^a_{m,\pm}=\Op(((|\xi '|^2+m^2)^\frac12 \pm i\xi _n)^a);  
$$
 where $\Xi ^a_{m,\pm} $ 
have exactly the same mapping properties as 
$\Xi ^a_\pm$,
 which is the case $m=1$. In particular, Theorem 3.1 holds
 with $\Xi ^a_\pm$ replaced by $\Xi ^a_{m,\pm} $.
It is seen as in \cite{G15}, Th.\ 4.2 and 4.4  that 
$$
r^+(-\Delta +m^2)u=r^+\Xi ^a_{m,-}e^+r^+ \Xi ^a_{m,+}u,
$$
when $u$ satisfies one of the mentioned hypotheses. Set $w=r^+ \Xi
^a_{m,+}u$,
 $w'= r^+\Xi ^a_{m,+}u'$.

We can then apply Theorem 3.1 to the integrals in the left-hand side
of (3.12), resp.\ the dualities in the left-hand side of (3.13), when
$u,u'$ satisfy the respective hypotheses there.
This gives e.g.\ under the weakest hypotheses (in $2^\circ$):
$$
\aligned
&\ang{r^+(-\Delta +m^2)^au,\partial_n 
u'}_{ \ol H^{\frac12-a+\varepsilon } ,\dot H^{a-\frac12-\varepsilon
}}
+\ang{\partial_nu,r^+(-\Delta +m^2)^au'}_{\dot H^{a-\frac12-\varepsilon
}, \ol H^{\frac12-a+\varepsilon } }\\
&=\ang{r^+\Xi _-^a e^+w,\partial_n 
u'}_{ \ol H^{\frac12-a+\varepsilon } ,\dot H^{a-\frac12-\varepsilon
}}+\ang{\partial_nu,r^+\Xi _-^a e^+w'}_{\dot H^{a-\frac12-\varepsilon}, \ol H^{\frac12-a+\varepsilon }} 
\\
&=2(\gamma _0w,\gamma _0w')_{L_2(\R^{n-1})}+\ang{w,\partial_n 
w'}_{ \ol H^{\frac12-\varepsilon } ,\dot H^{-\frac12+\varepsilon
}}+\ang{\partial_nw, 
w'}_{\dot H^{-\frac12+\varepsilon
},\ol H^{\frac12-\varepsilon }}.
\endaligned\tag3.14
$$
Let $w_k$ and $w'_k$ be sequences in $C^\infty _{(0)}(\crnp)=r^+C^\infty _0(\R^n)$ converging to $w$ resp.\ $w'$ in
$\ol H^{\frac12+\varepsilon }(\rnp)$ for $k\to\infty $; then $\gamma _0w_k\to \gamma _0w$
in $H^\varepsilon ({\Bbb
R}^{n-1})$ and $\partial_nw_k\to \partial_n w$ in $\ol H^{-\frac12+\varepsilon
}(\rnp)$, with
similar statements for $w'$. Now
$$
\aligned
\int_{\rnp}(w_k\partial_n \bar w'_k+\partial_n w_k\,\bar w'_k)\,
dx&=\int_{\rnp}\partial_n(w_k \bar w'_k)\, dx\\
&=-\int_{{\Bbb
R}^{n-1}}\gamma _0(w_k\bar w')\, dx'\to =-\int_{{\Bbb
R}^{n-1}}\gamma _0w\gamma _0\bar w'\, dx'. 
\endaligned
 $$
Thus the last two terms in (3.14) contribute with $-(\gamma _0w,\gamma _0w')_{L_2(\R^{n-1})}$, and we find that
$$
\aligned
&\ang{r^+\Xi _-^a e^+w,\partial_n 
u'}_{ \ol H^{\frac12-a+\varepsilon } ,\dot H^{a-\frac12-\varepsilon
}}+\ang{\partial_nu,r^+\Xi _-^a e^+w'}_{\dot H^{a-\frac12-\varepsilon}, \ol H^{\frac12-a+\varepsilon }}
\\
&=(\gamma _0w,\gamma _0w')_{L_2(\R^{n-1})}.
\endaligned
$$

Finally,
we recall from \cite{G15} Th.\ 5.1 that 
$$
\gamma _0w=\gamma _0(\Xi ^a _{m,+}u)=\gamma
_{a,0}u=\Gamma (a+1)\gamma _0(x_n^{-a}u).\tag3.15 
$$
Hence
$$
(\gamma _0w,\gamma _0w')_{L_2(\R^{n-1})}=\Gamma (a+1)^2\int_{{\Bbb
R}^{n-1}}\gamma _0(x_n^{-a}u)\,\gamma _0(x_n^{-a}\bar u')\, dx',
$$
and (3.13) follows. Under the hypotheses for $1^\circ$ and $3^\circ$
it can be written in the form (3.12).
\qed
\enddemo

The theorem also holds for $(-\Delta )^a$ itself (the case $m=0$), see Corollary 3.5 below.

We now turn to a 
general $\psi $do $P$
of order $2a$, elliptic avoiding a ray. The symbol is assumed to have the $a$-transmission property at the
 hyperplanes $\{x_n=c\}$, $c\in{\Bbb R}$ (as in (A.4) with $m=2a$, $\mu =a$):
$$
\partial_x^\beta \partial_\xi ^\alpha p_j(x,0,-\xi _n)=e^{\pi i(-j-|\alpha |)}\partial_x^\beta \partial_\xi ^\alpha p_j(x,0,\xi _n)\text{ for }j\in {\Bbb N}_0,\alpha ,\beta \in{\Bbb N}_0^n,\; |\xi |\ge 1;
$$
 this holds in particular when
 the symbol is even (cf.\ (2.40)). 

\proclaim{Theorem 3.4} Let $P$ be a classical $\psi $do on ${\Bbb R}^n$  of order $2a$ for some $0<a<1$, that is
elliptic avoiding a ray, with symbol having the $a$-transmission property at the
 hyperplanes $\{x_n=c\}$, $c\in{\Bbb R}$. Then the factorization index
 is $a$.

Let $s_0(x)=p_0(x,0,1)$, where $p_0$ is the principal
symbol, and let $P^{(n)}$ denote the commutator $[P,\partial_n]$;
$$
P^{(n)}=P\partial_n-\partial_nP;\text{ it has symbol }p^{(n)}=-\partial_{x_n}p,\tag3.16
$$
likewise  of order $2a$ and having the $a$-transmission property at the
 hyperplanes $\{x_n=c\}$.

For $u,u'\in H^{a(s)}(\crnp)$, $s\ge a+1$, there holds 
$$
\multline
\int_{\rnp} Pu\,\partial_n\bar
u'\,dx+\int_{\rnp}\partial_nu\,\overline{  P^*u'}\,dx\\=\Gamma (a+1)^2\int_{{\Bbb
R}^{n-1}}s_0\gamma _0(x_n^{-a}u)\,\gamma _0(x_n^{-a}\bar u')\,
dx'+\int_{\rnp}P^{(n)}u\, \bar u'\,dx.
\endmultline\tag3.17
$$

For $s\ge a+\frac12+\varepsilon $ (for some small $\varepsilon $), the formula holds with the integrals interpreted as
dualities:
$$
\multline
\ang{r^+ Pu,\partial_n
u'}_{\ol
H^{\frac12-a+\varepsilon }, \dot
H^{a-\frac12-\varepsilon}}+\ang{\partial_nu,{  P^*u'}}_{\dot H^{a-\frac12+\varepsilon},\ol
H^{\frac12-a-\varepsilon }}
\\=\Gamma (a+1)^2\int_{{\Bbb
R}^{n-1}}s_0\gamma _0(x_n^{-a}u)\,\gamma _0(x_n^{-a}\bar u')\,
dx'+\ang{r^+ P^{(n)}u,
u'}_{\ol
H^{\frac12-a+\varepsilon }, \dot H^{a-\frac12-\varepsilon }};
\endmultline\tag3.18
$$
the last term is a scalar product $(P^{(n)}u,u')_{L_2(\rnp)}$ when
$a\le \frac12$.  

In particular, when the symbol is independent of $x_n$, the term
with $P^{(n)}$ drops out.
\endproclaim

\demo{Proof} First let us account for the definition of the terms in
(3.17)--(3.18).
We
already have the information (3.8)--(3.9) on $u,u',\partial_nu$ and $\partial_nu'$.
If $u,u'\in H^{a(a+1)}(\crnp)$, we have the information (3.10).

By \cite{G15} Th.\ 4.2,
 $r^+P$ maps
 $H^{a(s)}(\crnp)$ continuously into $\ol H^{s-2a}(\rnp )$. When $s\ge
 a+1$, this
 is contained in $\ol H^{1-a}(\rnp )\subset L_2(\rnp )$, 
so $r^+Pu$ is an $L_2$-function. When $s\ge a+\frac12+\varepsilon $,
$r^+Pu\in \ol H^{\frac12-a+\varepsilon }(\rnp)$; in $L_2(\rnp)$ when $a\le\frac12$. 
The operator $P^{(n)}$, being of the same type as $P$, also has these
  mapping properties.

We see that for $s\ge a+1$, the first and last integrands in (3.17) are
functions.
For $s>a+\frac12$, the duality 
$$
\ang{r^+ Pu,\partial_n
u'}_{\ol
H^{\frac12-a+\varepsilon }, \dot H^{a-\frac12-\varepsilon }}.
$$
makes sense for small $\varepsilon $;
 here $r^+Pu$ is a function  when $a\le\frac12$, and $\partial_nu'$ is
 a function when $a>\frac12$. In the duality
$$
\ang{r^+ P^{(n)}u,
u'}_{\ol
H^{\frac12-a+\varepsilon }, \dot H^{a-\frac12-\varepsilon }},
$$
it is only $P^{(n)}u$ that may not be a function; it will be one when
 $a\le \frac12$.
(Observe also that since $a-\frac12\in \,]-\frac12,\frac12[\,$, $\ol
H^{\frac12-a+\varepsilon }\simeq \dot
H^{\frac12-a+\varepsilon }$  and $ \dot H^{a-\frac12-\varepsilon
 }\simeq  \ol H^{a-\frac12-\varepsilon }$ for small $\varepsilon >0$.)

The integral with $P^*$ is understood in a
similar way
(after conjugation).

In the right-hand sides of (3.17)--(3.18), the boundary values $\gamma
_0(x_n^{-a}u),\gamma _0(x_n^{-a}u')$ are defined as functions in
$H^{a+\frac12+\varepsilon -a-\frac12}(\R^{n-1})=H^{\varepsilon
}(\R^{n-1} )\subset L_2(\R^{n-1})$, by
\cite{G15}, Th.\ 5.1.

Now we turn to the proof of the formulas. The detailed arguments will
be given under the weakest regularity
hypothesis, namely $u,u'\in H^{a(a+\frac12+\varepsilon) }(\crnp)$.

In the reduction of the operators we  shall use  $\Lambda
^a_\pm$ (cf.\ (A.5)ff.)\ rather than $\Xi ^a_\pm$, in order to have
true $\psi $do's.
Then we write
$$
P=\Lambda _-^aQ\Lambda _+^a,\quad P^*=\Lambda _-^aQ^*\Lambda _+^a,\tag3.19
$$
where $Q=\Lambda _-^{-a}P\Lambda _+^{-a}$ is of order 0. It has the $0$-transmission
property at $\{x_n=0\}$, since $P$ is of type $a$, $\Lambda _-^{-a}$ is of type 0 and $\Lambda _+^{-a}$ is of type $-a$.
$Q$ is again elliptic avoiding a ray, since the symbols $\lambda _\pm^a$
of $\Lambda
_\pm^{a}$ are complex conjugates. An application
of Theorem 2.6 gives the factorization $q_0(x,\xi )=s_0(x)q_0^-(x,\xi
)q_0^+(x,\xi )$ with factors of order 0; then $p_0=\lambda _-^as_0q_0^-q_0^+\lambda _+^a$ with the plus-factor $q_0^+\lambda
_+^a$, and hence $P$ has factorization
index $a$. (We can normalize $\lambda ^a_\pm$ such that $\lambda
^a_\pm(0,1)=1$; then $s_0(x)=q_0(x,0,1)=p_0(x,0,1)$.)

Now construct $\psi $do's $Q_0 ^+$ and $Q_0 ^-$ from the symbols $q_0 ^+$
and $q_0 ^-$, 
and denote 
$$
Q-s_0 Q_0 ^-Q_0 ^+=  R_1, \quad  R=\Lambda
_-^a  R _1\Lambda _+^a. \tag3.20$$ 
Here $  R_1$ has 
order $-1  $, as a generalized $\psi $do, with symbol in $S^{-1 
}_{1,0}({\Bbb R}^{n},{\Bbb R}^{n-1},\Cal 
H_{-1})$. Indeed, $R_1$ has the symbol $q-q_0+q_0-(s_0q_0^-)\# q_0^+$,
where $q-q_0$ is a $\psi $do symbol of order $-1$ of type 0, and  
$$
q_0-(s_0q_0^-)\# q_0^+\sim s_0\sum _{|\alpha |\ge 1}\partial_\xi ^\alpha q_0^-D_{x}^\alpha q_0^+/\alpha !,
$$
where differentiation with respect to $\xi $ removes the term 1 in
$q_0^-$ and lowers the order, so that the resulting symbol is in 
 $S^{-1 }_{1,0}({\Bbb R}^{n},{\Bbb R}^{n-1},\Cal H_{-1})$.

For the main part of the operator $P_1=P-  R$ we use the factorization
$$
P_1=\Lambda _-^as_0 Q_0 ^-Q_0 ^+\Lambda _+^a=P^-P^+,\quad P^-=\Lambda _-^as_0 Q_0 ^-,\;
P^+=Q_0 ^+\Lambda _+^a; \tag3.21$$
here $P^-$ is a minus-operator, preserving support in $\crnm$,
and $P^+$ is a plus-operator, preserving support in $\crnp$.   Then we have the decompositions
$$
P=P^-{P^+}+  R,\quad P^*={P^+}^*{P^-}^*+  R^*.
$$

Let us first treat $P_1=P^-P^+$. We define
$$
\aligned
&w=r^+P^+u, \quad w'=r^+{P^-}^*u', \text{ then}\\
r^+P_1u&=r^+P^-e^+r^+P^+u=r^+P^-e^+w,\quad
r^+P^*u=r^+{P^+}^*e^+r^+{P^-}^*u'=r^+{P^+}^*e^+w',
\endaligned
$$
as in \cite{G15} Th.\ 4.2.
Here 
$r^+P_1u,r^+P_1^*u'\in \ol
H^{\frac12-a+\varepsilon }(\rnp)$, and,
 as noted further above,
$u,u'\in \dot H^{a+\frac12-\varepsilon }(\crnp)$ with
 $\partial_nu,\partial_nu'\in \dot
H^{a-\frac12-\varepsilon }(\crnp)$.

Define $v=r^+\Xi _+^au$, $v_1=r^+\Lambda  _+^au$, and recall that by the definition of
$H^{a(a+\frac12+\varepsilon )}(\crnp)$ in \cite{G15}, $$
v,v_1\in  \ol
H^{\frac12+\varepsilon }(\rnp),\text{ with }u=\Xi
_+^{-a}e^+v=\Lambda _+^{-a}e^+v_1.\tag3.22
$$ 
For $w$ we have that $
w=r^+Q_0^+\Lambda _+^au=r^+Q_0^+e^+v_1$.
Here $e^+v_1 \in e^+\ol H^{\frac12+\varepsilon }(\rnp)\subset \dot
H^{\frac12-\varepsilon '}(\crnp)$  (any $\varepsilon '>0$), which 
allows the conclusion that
$w\in \ol H^{\frac12-\varepsilon '}(\rnp)$, but we need to show that
$w\in \ol H^{\frac12+\varepsilon }(\rnp)$ (and similarly for $w'$). 
To do this, we shall use a (rough) parametrix $\widetilde
Q^-_0=\Op(1/q^-_0)$ of $Q^-_0$, cf.\
 Theorem 2.8. It is a minus-operator that satisfies 
$$
 \widetilde Q^-_0\, Q^-_0=I+R_2, \tag 3.23
$$
where $R_2$ is a minus-operator with symbol in $S^{-1
}_{1,0}({\Bbb R}^{n},{\Bbb R}^{n-1},\Cal H^-_{-1})$.
Denote
$r^+P_1u=f$, and recall that $f=r^+\Lambda _-^as_0Q_0^-e^+w$. Let
$$
w_1=r^+( \widetilde Q^-_0 s_0^{-1}\Lambda _-^{-a})e^+f;
$$
then since $r^+( \widetilde Q^-_0s_0^{-1} \Lambda _-^{-a})e^+\colon \ol
H^{s}(\rnp)\to \ol H^{s+a}(\rnp)$ for all $s$, $w_1\in \ol
H^{\frac12+\varepsilon }(\rnp)$. Now
$$
w_1-w=r^+( \widetilde Q^-_0s_0^{-1} \Lambda _-^{-a})e^+r^+(  \Lambda
_-^as_0Q_0^-)e^+w-w
=r^+( \widetilde Q^-_0 Q_0^-)e^+w-w=r^+R_2e^+w;
$$
where we used that $e^+r^+$ in the middle can be left out since the
operators are minus-operators, that $\Lambda ^{-a}_-\Lambda ^a_-=I$,
and that (3.23) holds. Here $r^+R_2e^+$ maps $w$ into $\ol
H^{\frac32-\varepsilon ' }(\rnp)$, by Theorem 2.4. It follows that 
$w\in \ol H^{\frac12+\varepsilon }(\rnp)$. A similar proof shows this for  $w'$.

Now we can write
$$
I_1\equiv \ang{r^+P_1u,\partial_n u'}_{\ol
H^{\frac12-a+\varepsilon }, \dot H^{a-\frac12-\varepsilon}}
=\ang{r^+P^-e^+w,\partial_n u'}_{\ol H^{\frac12-a+\varepsilon
}, \dot H^{a-\frac12-\varepsilon }}.
$$
Since $r^+P^-e^+\colon \ol H^{\frac12+\varepsilon }(\rnp)\to\ol H^{\frac12-a+\varepsilon }(\rnp)$ and
${P^-}^*\colon \dot H^{a-\frac12-\varepsilon }(\crnp)\to \dot H^{-\frac12-\varepsilon }(\crnp)$ are adjoints,
$$
I_1=\ang {w, {P^-}^*\partial_nu'}_{\ol H^{\frac12+\varepsilon }, \dot H^{-\frac12-\varepsilon }}.
$$
We use here that $u'$ is zero at $x_n=0$, so that
$\partial_nu'=e^+r^+\partial_nu'$ (one may identify
$\partial_nu'$ with $r^+\partial_nu'$).

The distribution $ {P^-}^*\partial_nu'\in \dot H^{-\frac12-\varepsilon }(\crnp)$ is
rewritten as follows:
$$
{P^-}^*\partial_nu'=\partial_n{P^-}^*u'+[{P^-}^*,\partial_n]u'=\partial_ne^+w'+{{P^-}^*}^{(n)}u',
$$
with the notation $[{P^-}^*,\partial_n]={{P^-}^*}^{(n)}$ as in (3.16).
Here, as in Theorem 3.1, 
$$
\partial_ne^+w'=\gamma _0(w')\otimes \delta (x_n)+e^+\partial_nw',
$$
where we moreover note that since $w'\in \ol H^{\frac12+\varepsilon }(\crnp)$, 
$e^+\partial_nw'$ is not just in $\dot
H^{-\frac12-\varepsilon }(\crnp)$,
but is in  $\dot
H^{-\frac12+\varepsilon }(\crnp)\simeq\ol
H^{-\frac12+\varepsilon }(\rnp)
$, and $\gamma _0w'\in H^{\varepsilon }(\R^{n-1})$.
Insertion of the expressions in $I_1$ and integration by parts as in
Theorem 3.1 gives:
$$
\aligned
I_1&=\ang {w,\gamma _0(w')\otimes \delta
(x_n)+e^+\partial_nw'+{{P^-}^*}^{(n)}u'}
_{\ol H^{\frac12-\varepsilon }, \dot H^{-\frac12+\varepsilon }}\\
&=(\gamma _0w,\gamma
_0w')_{L_2(\R^{n-1})}+\ang{w,\partial_nw'}_{\ol H^{\frac12-\varepsilon
}, \dot H^{-\frac12+\varepsilon }}
+\ang{w,{{P^-}^*}^{(n)}u'}
_{\ol H^{\frac12-\varepsilon }, \dot H^{-\frac12+\varepsilon }}.
\endaligned\tag3.24
$$

It is shown in the same way (in fact it can be concluded from the
above by interchanging $P_1$ and $P_1^*$, $u$ and $u'$, and conjugating), that
$$
I_2\equiv \ang{\partial_nu, r^+P_1^*u'}_{\dot H^{a-\frac12+\varepsilon},\ol
H^{\frac12-a-\varepsilon }}, 
$$
satisfies
$$
I_2=(\gamma _0w,\gamma
_0w')_{L_2(\R^{n-1})}+\ang{\partial_nw,w'}_{\dot
H^{-\frac12+\varepsilon },\ol H^{\frac12-\varepsilon}} 
+\ang{{P^+}^{(n)}
u,w'}_{\dot H^{-\frac12+\varepsilon },\ol H^{\frac12-\varepsilon}},\tag3.25 
$$
where  ${P^+}^{(n)}$  stands for $[P^+,\partial_n]$ as in (3.16).

Taking the two contributions together, we find that
$$
\aligned
I_1+I_2&=2(\gamma _0w,\gamma
_0w')_{L_2(\R^{n-1})}+\ang{w,\partial_nw'}_{\ol H^{\frac12-\varepsilon
}, \dot H^{-\frac12+\varepsilon }}+\ang{\partial_nw,w'}_{\dot
H^{-\frac12+\varepsilon },\ol H^{\frac12-\varepsilon}} \\
&\quad +\ang{w,{{P^-}^*}^{(n)}u'}_{\ol H^{\frac12-\varepsilon }, \dot H^{-\frac12+\varepsilon }}+\ang{{P^+}^{(n)}
u,w'}_{\dot H^{-\frac12+\varepsilon },\ol H^{\frac12-\varepsilon}}\\
&=(\gamma _0w,\gamma
_0w')_{L_2(\R^{n-1})}+I_3,\text{ where}\\
I_3&=\ang{P^+u,{{P^-}^*}^{(n)}u'}_{\ol H^{\frac12-\varepsilon} , \dot H^{-\frac12+\varepsilon }}+\ang{{P^+}^{(n)}
u,{P^-}^*u'}_{\dot H^{-\frac12+\varepsilon },\ol H^{\frac12-\varepsilon}};
\endaligned\tag 3.26
$$
here we used the calculation after (3.14) to reduce the first line to
a single boundary integral, and collected the last two
terms in $I_3$. This will now be further reduced.

Observe that ${P^+}^{(n)}$ has symbol equal to $-\partial_{x_n}$ of the symbol
of $P^+=\Lambda _+^aQ_0^+$, so it is a plus-operator, continuous from
$\dot H^s(\crnp)$ to $\dot H^{s-a}(\crnp)$ for all $s$, with an
adjoint  $r^+({P^+}^{(n)})^*e^+$ going from  $\ol H^{a-s}(\rnp)$ to $\ol H^{-s}(\rnp)$ for all
$s$.
${{P^-}^*}^{(n)}$ has similar properties. In particular,
${P^+}^{(n)}u={P^+}^{(n)}\Xi _+^{-a}e^+v$
is in 
$ \dot
H^{\frac12-\varepsilon' }(\crnp)$, cf.\ (3.22), and so is
${{P^-}^*}^{(n)}u'$,
 so the dualities in $I_3$ identify with
$L_2(\rnp)$-scalar products:
$$
I_3=(P^+u,{{P^-}^*}^{(n)}u')_{L_2(\rnp)}
+({P^+}^{(n)}u,{P^-}^*u')_{L_2(\rnp)}.
$$
The adjoint of ${{P^-}^*}^{(n)}$ is
 $r^+{P^-}^{(n)}e^+$, since
 $[{P^-}^*,\partial_n]^*=[P^-,\partial_n]$. Then in view of the
 mapping properties,
$$
\aligned
I_3&=\ang{r^+{P^-}^{(n)}e^+r^+P^+u,u'}_{\ol H^{\frac12-a+\varepsilon} , \dot H^{a-\frac12-\varepsilon }}+\ang{r^+P^-e^+r^+{P^+}^{(n)}
u,u'}_{\ol H^{\frac12-a+\varepsilon },\dot H^{a-\frac12-\varepsilon}}\\
&=\ang{r^+({P^-}^{(n)}e^+r^+P^++r^+P^-e^+r^+{P^+}^{(n)})
u,u'}_{\ol H^{\frac12-a+\varepsilon },\dot H^{a-\frac12-\varepsilon}}.
\endaligned\tag3.27
$$

We now use moreover, that $$
r^+P^-e^+r^+{P^+}^{(n)} u=r^+P^-{P^+}^{(n)}
u,\quad r^+{P^-}^{(n)}e^+r^+P^+u=
r^+{P^-}^{(n)}P^+u
$$ 
(because of the support-preserving properties, as
in \cite{G15} Th.\ 4.2), so that
$$
I_3=\ang{r^+(P^-{P^+}^{(n)}+{P^-}^{(n)}P^+)u, u'}_{\ol H^{\frac12-a+\varepsilon },\dot H^{a-\frac12-\varepsilon}}.
$$
 Here we can perform a little calculation on the $\psi $do's on $\R^n$:
$$
\aligned
P^-{P^+}^{(n)}+{P^-}^{(n)}P^+&=P^-P^+\partial_n-P^-\partial_nP^++P^-\partial_nP^+
-\partial_nP^-P^+\\
&=P^-P^+\partial_n-\partial_nP^+P^-=P_1^{(n)},
\endaligned\tag3.28
$$
showing that in fact
$$
I_3=\ang{r^+P_1^{(n)}u,u'}_{\ol H^{\frac12-a+\varepsilon },\dot H^{a-\frac12-\varepsilon}}.
$$
Inserting this in (3.26), we reach the conclusion that
$$
I_1+I_2=(\gamma _0(P^+u),\gamma _0({P^-}^*u'))_{L_2({\Bbb R}^{n-1})}+\ang{r^+P_1^{(n)}u,u'}_{\ol H^{\frac12-a+\varepsilon },\dot H^{a-\frac12-\varepsilon}}.\tag3.29
$$

The boundary term can be further clarified as follows:
Let $v=r^+\Xi _+^au$ as in (3.22). We know from \cite{G15} (cf.\ e.g.\ Cor.\ 5.3) that $\gamma _0v=\Gamma
(a+1)\gamma _0((x_n^{-a})u)$. In view of Theorem 2.6, we have that
$$
q_0^+ (x,\xi ',\xi _n)=1+f(x,\xi ',\xi _n)
$$
where $f$ is in $\Cal H^+$ as a function of $\xi _n$; 
$f\in S^0_{1,0}({\Bbb R}^n,{\Bbb R}^{n-1},\Cal H^+)$. Hence $$
Q_0^+=I+F,\quad F=\Op(f).
$$
Moreover, $\Lambda
_+^a=(1+\Psi )\Xi _+^a$, where $\Psi $ has symbol $\psi (\xi )$ in $\Cal H^+$ with respect to $\xi _n$, cf.\
\cite{G15} (1.16) and Lemma 6.6; $\psi \in S^0_{1,0}({\Bbb R}^n,{\Bbb
R}^{n-1},\Cal H^+)$. It follows that
$$
Q_0^+\Lambda _+^a=(I+F)(I+\Psi )\Xi ^a_+=(I+F_1)\Xi _+^a,
$$
where $F_1$ has symbol $f_1\in S^0_{1,0}({\Bbb R}^n,{\Bbb
R}^{n-1},\Cal H^+)$. (One could also deal with the factors $I+F$ and
$I+\Psi $ in two successive steps, to avoid using Leibniz products.) 
By the rules of the Boutet de Monvel calculus,
$$
\gamma _0(F_1v)=(2\pi )^{-n}\int_{{\Bbb R}^{n-1}}e^{ix'\cdot \xi
'}\int_{{\Bbb R}}f_1(x,\xi ',\xi _n)\Cal F(e^+v(x',x_n))\,d\xi _nd\xi '=0.\tag3.30
$$
(Briefly recalled, the reason is that both $f_1$ and $\Cal F(e^+v(x',x_n))$ are in $\Cal H^+$ as
functions of $\xi _n$ --- the latter because $e^+v$ is supported in
$\crnp$; then their product is $O(\ang{\xi _n}^{-2})$ and holomorphic
in ${\Bbb C}_-$, so the integral over ${\Bbb R}$ can be transformed to
a closed contour in ${\Bbb C}_-$ and therefore vanishes.)
It follows that
$$
\gamma _0(P^+u)=\gamma _0(Q_0^+\Lambda _+^au)=\gamma _0(v+F_1v)=\gamma _0v
=\gamma _0(\Xi _+^au)=\Gamma (a+1)\gamma _0(x_n^{-a}u). \tag3.31$$

As a slight variant, we also have, with $v'=\Xi _+^au'$:
$$
\aligned
\gamma _0({P^-}^*u')&=\gamma _0({Q_0^-}^*\bar s_0\Lambda _+^au')=\gamma
_0({Q_0^-}^*\bar s_0(I+\Psi )v')=\gamma _0(\bar s_0v'+F_2v')\\
&=\gamma _0(\bar s_0v')=\bar s_0\Gamma (a+1)\gamma _0(x_n^{-a}u'),
\endaligned\tag3.32
$$
where ${Q_0^-}^*\bar s_0(I+\Psi )=\bar s_0I+F_2$, and also  $F_2$ has symbol in $\Cal H^+$ w.r.t.\ $\xi _n$, hence
does not contribute. (Recall that $s_0$ is a function of $x$, namely
$s_0(x)=q_0(x,0,1)=p_0(x,0,1)$; in the final formula it is just its
value on $\{x_n=0\}$ that enters.)

We conclude that
$$
(\gamma _0(P^+u),\gamma _0({P^-}^*u'))_{L_2(\R^{n-1})}
=\Gamma (a+1)^2\int_{{\Bbb R}^{n-1}}\gamma _0(x_n^{-a} u)\,s_0\gamma
_0(x_n^{-a}\bar u')\,dx',\tag3.33
$$
whereby
$$
I_1+I_2=\Gamma (a+1)^2\int_{{\Bbb R}^{n-1}}s_0\,\gamma _0(x_n^{-a} u)\,\gamma
_0(x_n^{-a}\bar u')\,dx'+\ang{r^+P_1^{(n)}u,u'}_{\ol H^{\frac12-a+\varepsilon },\dot H^{a-\frac12-\varepsilon}}.
\tag3.34
$$


Finally, we must also treat the contribution from $  R=\Lambda
_-^a  R _1\Lambda _+^a$. 
As already noted, the symbol $r_1(x,\xi )$ of $R_1$ is in
 $ \Cal H_{-1}$ as a function of $\xi _n$, so we can apply the
projections $h^+$ and $h^-$, decomposing
$$
r_1(x,\xi )=r_1^+(x,\xi )+ r_1^-(x,\xi ),\quad r_1^\pm\in  S^{-1
}_{1,0}({\Bbb R}^{n},{\Bbb R}^{n-1},\Cal 
H_{-1}^\pm).\tag3.35
$$
 Denote the hereby defined operators $  R_1^\pm$; $  R_1= 
 R_1^++  R_1^-$. Then when we set $S^-=\Lambda _-^a  R ^-_1,\,   S^+=
 R^+ _1\Lambda _+^a
$, $  R$ is decomposed as
$$
R=\Lambda
_-^a  R ^-_1\Lambda _+^a+\Lambda
_-^a  R^+ _1\Lambda _+^a=   S^-\Lambda _+^a+\Lambda
_-^a  S^+ ;\tag3.36 
$$
a sum of two operators that are products of a minus-operator and a
plus-operator. To each of these products, we can apply
the same method as we did to $P^-P^+$. This reduces the corresponding
integrals to scalar products over the boundary plus commutator contributions:
$$
\aligned
I_4&\equiv \ang{r^+Ru,\partial_n u'}_{\ol
H^{\frac12-a+\varepsilon }, \dot H^{a-\frac12-\varepsilon}}+\ang{\partial_nu, r^+R^*u'}_{\dot H^{a-\frac12+\varepsilon},\ol
H^{\frac12-a-\varepsilon }}\\
&=(\gamma _0(\Lambda _+^au),\gamma _0({{  S^-}^*u'}))_{L_2(\R^{n-1})}+\ang{\Lambda _+^au,{{S^-}^*}^{(n)}u'}_{\ol H^{\frac12-\varepsilon} , \dot
H^{-\frac12+\varepsilon }}\\
&\quad+(\gamma _0(  S^+u),\gamma _0({\Lambda
 _+^au'}))_{L_2(\R^{n-1})}+\ang{{S^+}^{(n)}
u,\Lambda _+^au'}_{\dot H^{-\frac12+\varepsilon },\ol
H^{\frac12-\varepsilon}}
.\endaligned
\tag3.37
$$
(The dualities in the second and third line reduce to $L_2$-scalar products since $S^-$ and
$S^+$ are of negative order.)
Since ${ R_1^-}^*$ and $  R_1^+$ have symbols in $\Cal  H^+$ as
functions of $\xi _n$, the boundary values of ${S^-}^*u'$ and $S^+u$ are zero, so only
the commutator terms survive. These are reduced in a similar way as in
the treatment of $P_1$, to give
$$
I_4=(r^+R^{(n)}u,u')_{L_2(\rnp)}.
$$

Collecting all the terms, we find (3.18). As accounted for in the
beginning of the proof, it can be written in the
form (3.17) when $u,u'\in H^{a(a+1)}(\crnp)$.
\qed
\enddemo

One can in particular conclude:

\proclaim{Corollary 3.5} For $u,u'\in H^{a(s)}(\crnp)$ with $s>a+\frac12$,
$$
\multline
\ang{r^+(-\Delta)^au,\partial_n 
u'}_{ \ol H^{\frac12-a+\varepsilon } ,\dot H^{a-\frac12-\varepsilon
}}
+\ang{\partial_nu,r^+(-\Delta )^au'}_{\dot H^{a-\frac12-\varepsilon
}, \ol H^{\frac12-a+\varepsilon } }\\
=\Gamma (a+1)^2\int_{{\Bbb
R}^{n-1}}\gamma _0(x_n^{-a}u)\,\gamma _0(x_n^{-a}\bar u')\, dx'.
\endmultline\tag3.38
$$
for small $\varepsilon >0$. When $s\ge a+1$, this can be written as
$$
\multline
\int _{\rnp}(-\Delta )^au\,\partial_n\bar u'\, dx+\int
_{\rnp}\partial_nu\,{(-\Delta )^a\bar u'}\, dx
=\Gamma (a+1)^2\int_{{\Bbb
R}^{n-1}}\gamma _0(x_n^{-a}u)\,\gamma _0(x_n^{-a}\bar u')\, dx'.
\endmultline
\tag 3.39
$$

\endproclaim

\demo{Proof} Write
$$
(-\Delta )^a=P+\Cal S,\text{ where }P=\Op(\eta (\xi )|\xi |^{2a}),\;
\Cal S=\Op((1-\eta (\xi ))|\xi |^{2a});\tag3.40
$$
$\eta (\xi )$ denoting an excision function as in (2.11). Then $P$
satisfies the hypotheses of Theorem 3.4, so (3.38) holds for this 
operator.

Now consider $\Cal S$. Its symbol $s(\xi )=(1-\eta (\xi ))|\xi
|^{2a}$ is bounded and supported in $\ol B_1=\{|\xi |\le 1\}$. The same
holds for all the symbols $s_\alpha =\xi ^\alpha (1-\eta (\xi ))|\xi
|^{2a}$, $\alpha \in{\Bbb N}_0^n$, so they all define bounded
operators in $H^t(\R^n)$,  for all
$t\in \R$. Since $\Op (s_\alpha )=D^\alpha \Cal S=\Cal SD^\alpha
$,
we see that $\Cal S$  and
its compositions with $D^\alpha $ are smoothing operators, going from 
$H^\infty (\R^n)=\bigcup_t H^t(\R^n)$  to $H^{-\infty }(\R^n)=\bigcap_t H^t(\R^n)$.

Recall from (3.9) that $u\in \dot H^{a+\frac12-\varepsilon '}(\crnp)$,
$\partial_nu\in \dot H^{a-\frac12-\varepsilon '}(\crnp)$ for any
$\varepsilon '>0$; here we can choose $\varepsilon '$ so that
$\sigma =a-\frac12-\varepsilon '\in \,]-\frac12,\frac12[\,$.
Then $\Cal Su\in H^{-\infty }(\R^n)$; and 
$$
\ang{r^+\Cal Su, \partial_nu'}_{ \ol H^{-\sigma  } ,\dot H^{\sigma }}=\ang{\Cal Su, \partial_nu'}_{  H^{-\sigma  }(\R^n) , H^{\sigma }(\R^n)},\tag3.41
$$
since $\Cal Su=e^+r^+\Cal Su+e^-r^-\Cal Su$, where the terms are in
$\dot H^{|\sigma |}=\ol H^{\,|\sigma |} $ over $\crnp$ resp.\ $\crnm$,
and $e^-r^-\Cal Su$ vanishes on $\partial_nu'$. In the last expression in
(3.41), $\partial_n$ can be moved to the left-hand side with a minus,
and $\Cal S$ can be moved to the right-hand side replaced by $\Cal
S^*$ (of the same type), with suitable adaptation of the duality indications. 
Then we find that 
$$
\ang{r^+\Cal Su, \partial_nu'}_{ \ol H^{-\sigma  } ,\dot H^{\sigma }}+\ang{\partial_nu,r^+\Cal S^*u'}_{ \dot H^{\sigma  } ,\ol H^{-\sigma }}=0,
$$
and when this is added to
the integration by parts formula for $P$, we find (3.38).

When $s\ge a+1$, then $u$ and $\partial_nu$ are functions, and so are $r^+Pu$
and $\Cal Su$, with similar statements for $u'$. Then the formula can
be written as in (3.39).
\qed

\enddemo

\head 4. Integration by parts over bounded smooth domains \endhead

In this part, we consider a classical $\psi $do $P$ of order $2a$ on
${\Bbb R}^n$ and its restriction to a bounded smooth subset $\Omega
$. Assuming that the symbol is even (cf.\ (2.40)), we have that it
satisfies the $a$-transmission condition in any direction at all
points, hence at the boundary of
any choice of $\Omega $. The indications $r^\pm$ and $e^\pm$ now pertain
to the embedding $\Omega \subset \R^n$.

We begin with a simple integration-by-parts formula, that can be shown
 by reduction to operators of order 0.

\proclaim{Theorem 4.1} Let $P$ be a classical $\psi $do on ${\Bbb R}^n$ of order
$2a$ for some $a>0$, with even symbol. Then for $u,u'\in
H^{a(s)}(\comega)$, $s\ge a$,
$$
\ang{ r^+ Pu,   u'}_{\ol H^{-a}(\Omega ),\dot H^{a}(\comega)} -
\ang{u,{r^+ P^*u'}}_{\dot H^a(\comega),\ol H^{-a}(\Omega )}=0.\tag4.1
$$
When  $s\ge 2a$, this can also be written
$$
\int_{\Omega}  Pu\, \bar u'\,dx - \int_{\Omega }u\,\overline{ P^*u'}\,dx=0.\tag4.2
$$

\endproclaim

\demo{Proof} We shall apply the families of order-reducing operators $\Lambda
_+^{(t)}$ and $\Lambda _{-,+}^{(t)}$, $t\in{\Bbb R}$, introduced in
\cite{G15} and recalled in the Appendix, chosen such that $\Lambda _{-,+}^{(t)}\colon
\ol  H_{p}^s(\Omega  )\to
\ol  H_{p}^{s-a}(\Omega  )$ and
$\Lambda _+^{(t)}\colon \dot H_{p'}^{a-s}(\comega )\to \dot
 H_{p'}^{-s}(\comega )$
 are adjoints for all $s\in{\Bbb R}$. 
Recall that $H_p^{a(s)}(\comega)=\Lambda
 _+^{(-a)}e^+\ol  H_p^{s-a}(\Omega )$. We restrict the attention to the case $p=2$.

Since $P$ is even, it has the $a$-transmission property at any
boundary; then the operator 
$$
Q=\Lambda _-^{(-a)}P\Lambda _+^{(-a)},\tag4.3
$$
is  a $\psi $do of order 0 having the $0$-transmission property at
the boundary of $\Omega $. Recall that $r^+P$ maps
 $H^{a(s)}(\comega)$ continuously into $
 \ol  H^{s-2a}(\Omega )$ for all $s> a-\frac12$, cf.\ \cite{G15} Th.\ 4.2.
 
Let
$$
w=r^+\Lambda _+^{(a)}u, \quad w'=r^+\Lambda _+^{(a)}u ';
$$
they are in $
\ol H^{s-a}(\Omega )$, which identifies with a subset of $ L_2(\Omega
)$ since $s\ge a$.
Then
$$
u=\Lambda _+^{(-a)}e^+w,\, u'=\Lambda _+^{(-a)}e^+w' \in \dot H^a(\comega)
$$
(using that $\Lambda ^{(-a)} _+$ lifts $e^+L_2(\Omega )$ to $\dot H^a(\comega )$),
and $r^+Pu$, $r^+P^*u'\in \ol H^{s-2a}(\Omega )\subset \ol
H^{\,-a}(\Omega )$.
(Since $u$ is an $L_2$-function supported in $\comega$, we identify it
with $r^+u$.)
Moreover (cf.\ \cite{G15}), 
$$
\aligned
r^+Pu&=r^+\Lambda _{-}^{(a)}e^+r^+Q\Lambda _+^{(a)}u=r^+\Lambda _{-}^{(a)}e^+r^+Qe^+w,\\
r^+P^*u'&=r^+\Lambda _{-}^{(a)}e^+r^+Q^*\Lambda _+^{(a)}u'=r^+\Lambda _{-}^{(a)}e^+r^+Q^*e^+w'.
\endaligned
\tag4.4
$$

Now since $r^+\Lambda ^{(a)}_-e^+$ and $\Lambda ^{(a)}_+$ are adjoints,
$$
\ang{r^+Pu,  u'}_{\ol H^{-a}, \dot H^a}=\ang{r^+\Lambda
_{-}^{(a)}e^+r^+Qe^+w,{\Lambda _+^{(-a)}w'}}_{\ol H^{-a}, \dot
H^a}=(r^+Qe^+w,w')_{L_2(\Omega )}.
$$
There is a similar formula for $P^*$, so we find
$$
\aligned
&\ang{r^+Pu,  u'}_{\ol H^{-a},\dot H^a}-\ang{u,{ r^+P^*u'}}_{\dot H^a,
\ol H^{-a}}=(r^+Qe^+w,w')_{L_2(\Omega )}-(w,r^+Q^*e^+w')_{L_2(\Omega )}.
\endaligned\tag4.5
$$
Since $Q$ is of order 0, the adjoint of $r^+Qe^+$ in $L_2(\Omega )$ is
$r^+Q^*e^+$, and
$$
(r^+Qe^+w,w')_{L_2(\Omega )}-(w,r^+Q^*e^+w')_{L_2(\Omega )}=0.\tag4.6
$$
This shows (4.1). 

When $s\ge 2a$, $r^+Pu$ and $r^+P^*u'\in L_2(\Omega )$, so
the formula can be written as in (4.2).\qed 
\enddemo

The formula can be extended to suitable  $L_p,
L_{p'}$-dualities.

Our main aim is to show extensions of the integration-by-parts formula in Theorem 3.4 to the curved
situation.

First there is a result in the spirit of Theorem 3.1.

\proclaim{Theorem 4.2} Let $P^-$ be an operator of order $a$ (i.e.,
continuous from $H^s(\R^n)$ to $H^{s-a}(\R^n)$ for all $s\in{\Bbb R}$)
such
that $r^+P^-e^+$ maps $ \ol
H^{s }(\Omega )$ continuously to $\ol
H^{s-a }(\Omega )$ with 
adjoint ${P^-}^*\colon \dot H^{a-s}(\comega)\to \dot
H^{-a}(\comega)$ for all $s\in{\Bbb R}$. Assume that the commutator 
$$
{P^-}^{(j)}=P^-\partial_j-\partial_jP^-
$$
has similar mapping properties.
 Let
$w, w'\in \ol H^s(\comega)$ with $s\ge\frac12+\varepsilon $ for some small
$\varepsilon >0$, and assume that $w'=r^+{P^-}^*u'$ for some $u'\in
H^{a(s+a)}(\comega)$ with ${P^-}^*u'=e^+w'$. Then
$$
\multline
\ang{r^+P^-e^+w,\partial_ju'}_{\ol H^{\frac12-a+\varepsilon },\dot
H^{a-\frac12-\varepsilon}}
=
\int_{\partial\Omega }\nu _j\,\gamma _0w\,\gamma _0\bar w'\, d\sigma
+
\ang{w,\partial_jw'}_{\ol H^{\frac12-\varepsilon },\dot
H^{-\frac12+\varepsilon}}\\
+(w, {{P^-}^*}^{(j)}u')_{L_2(\Omega )},
\endmultline\tag4.7
$$
where $\nu _j(x)$ is the $j$'th component of the interior normal
vector $\nu (x)$ at $x\in\partial\Omega $.
\endproclaim

\demo{Proof} Recall the standard  Gauss-Green formula
$$
-\int_{\Omega }\partial_j\varphi \, dx=\int_{\partial\Omega }\nu _j\gamma
 _0\varphi \, d\sigma ,\tag 4.8
$$
where $\gamma _0\varphi $ is the restriction of $\varphi $ to
$\partial\Omega $ and $d\sigma $ is the induced measure on
$\partial\Omega $; it holds
for sufficiently regular functions $\varphi $. We can write it as a
distribution formula on $\R^n$ (with sesquilinear duality):
$$
\ang{\partial_j1_{\Omega }, \varphi }_{\R^n}=-\ang{1_\Omega
,\partial_j\varphi }_{\R^n}=\ang {1,\nu _j\widetilde \gamma _0\varphi
}_{\partial\Omega }\text{ for }\varphi \in C_0^\infty (\R^n),\tag4.9
$$
 where the last brackets is a duality over $\partial\Omega $
 consistent with the scalar product in $L_2(\partial\Omega ,d\sigma
 )$.
For accuracy, we denote by $\widetilde\gamma _0$ the restriction
 operator going from functions on $\R^n$ to functions on
 $\partial\Omega $ (sometimes called the two-sided trace operator); it
 is this one that has nice adjoint properties. In fact, 
$$
\widetilde \gamma _0\colon H^{s }({\Bbb R}^n)\to H^{s-\frac12
}(\partial\Omega )\text{ has an adjoint }
\widetilde \gamma _0^*\colon H^{\frac12-s }(\partial\Omega )\to H^{-s }({\Bbb
R}^{n})\text{ for }s>\tfrac12, \tag4.10
$$
and (4.9) shows that $\partial_j1_{\Omega }=\widetilde \gamma _0^*\nu
_j$.

There is also a version with two functions $W$ and $\varphi $:
When $W
\in C_0^\infty (\R^n)$, $\partial_j (1_\Omega W)=(\partial_j1_{\Omega
})W+1_{\Omega }\partial_jW$, so
$$
\multline
\ang{\partial_j(1_{\Omega }W)-1_\Omega \partial_jW,\varphi }_{\R^n}=
\ang{(\partial_j1_{\Omega })W,\varphi }_{\R^n}=\ang{\partial_j1_{\Omega
},\ol W\varphi }_{\R^n}=\ang{1,\nu _j\widetilde\gamma _0(\ol W\varphi
)}_{\partial\Omega }\\
=\ang{1,\nu
_j\widetilde\gamma _0(\ol W)\,\widetilde \gamma_0(\varphi )}_{\partial\Omega }
=\ang{\nu _j\widetilde\gamma _0( W),\widetilde\gamma _0\varphi }_{\partial\Omega }=\ang{\widetilde \gamma _0^*(\nu _j\widetilde\gamma _0( W)),\varphi }_{\R^n},
\endmultline
$$
showing that 
$$
\partial_j(1_{\Omega }W)=1_\Omega \partial_jW +\widetilde \gamma _0^*(\nu _j\widetilde\gamma _0( W)).
$$
Setting $r^+W=w$, we find the formula
$$
\partial_je^+w=e^+ \partial_jw +\widetilde \gamma _0^*(\nu _j\gamma _0
w).\tag 4.11
$$
It extends by continuity to more general functions, namely $w\in \ol
H^{\frac12+\varepsilon }(\Omega )$ with $\gamma _0w\in H^{\varepsilon
}(\partial\Omega )$.

For the left-hand side in (4.7)
 we then find:
$$
\aligned
&\ang{r^+P^-e^+w,\partial_ju'}_{\ol H^{\frac12-a+\varepsilon },\dot
H^{a-\frac12-\varepsilon}}=\ang{w,{P^-}^*\partial_ju'}_{\ol H^{\frac12+\varepsilon },\dot
H^{-\frac12-\varepsilon}}\\
&=\ang{w,\partial_j{P^-}^*u'+{{P^-}^*}^{(j)}u'}_{\ol H^{\frac12+\varepsilon },\dot
H^{-\frac12-\varepsilon}}=\ang{w,\partial_je^+w'+{{P^-}^*}^{(j)}u'}_{\ol H^{\frac12+\varepsilon },\dot
H^{-\frac12-\varepsilon}}\\
&=\ang{w,e^+ \partial_jw' +\widetilde \gamma _0^*(\nu _j\gamma _0
w')+{{P^-}^*}^{(j)}u'}_{\ol H^{\frac12+\varepsilon },\dot
H^{-\frac12-\varepsilon}}\\
&=\int_{\partial\Omega }\nu _j\gamma _0w\gamma _0
\bar w'\,d\sigma +
\ang{w,e^+ \partial_jw'}_{\ol H^{\frac12-\varepsilon },\dot
H^{-\frac12+\varepsilon}}+(w, {{P^-}^*}^{(j)}u')_{L_2(\Omega )}.
\endaligned
$$
Here we used the information on adjoints and inserted (4.11) applied
to $w'$; the dualitiy indications could be changed since $e^+\partial_jw' $
and ${{P^-}^*}^{(j)}u'$ lie in better spaces $\dot
H^{-\frac12+\varepsilon }$, resp\. $\dot
H^{\frac12+\varepsilon }$.
\qed
\enddemo

To treat the full problem, we shall use local coordinates.

Let $\Omega $ be a smooth bounded subset of ${\Bbb R}^n$. Then 
 $\comega$ has a finite cover by bounded open sets $U_0,\dots, U_{I_0}$ with
 diffeomorphisms
$\kappa _i\colon U_i\to V_i$, $V_i$ bounded open in ${\Bbb R}^n$, such
 that $U_i\cap \Omega $ is mapped to $V_i\cap\rnp$ and
 $U_i\cap\partial\Omega $ is mapped to $V_i\cap\partial\crnp$; as
 usual we write $\partial\crnp={\Bbb R}^{n-1}$. When $P$ is a $\psi
 $do on ${\Bbb R}^n$, its application to functions
 supported in $U_i$ carries over to functions on $V_i$ as a $\psi $do
 $\widetilde P^{(i)}$ defined by 
$$
\widetilde P^{(i)}v=P(v\circ \kappa _i)\circ \kappa _i^{-1},\quad v\in C_0^\infty (V_i).\tag4.12
$$

\example{Remark 4.3} A useful choice near $\partial\Omega $ is where we provide the $(n-1)$-dimensional manifold $\partial\Omega $ with coordinate
charts $\kappa '_i\colon U'_i\to V'_i\subset {\Bbb R}^{n-1}$,
$i=1,\dots, I_0$, and consider a tubular
neighborhood $\Sigma _r=\{x'+t\nu (x')\mid x'\in\partial\Omega , |t|<r\}$,
where $\nu (x')=(\nu _1(x'),\dots,\nu _n(x'))$ is the interior normal
to $\partial\Omega $ at $x'\in\partial\Omega $, and $r$ is taken so
small that the mapping $ x'+t\nu (x')\mapsto (x',t)$ is a
diffeomorphism from $\Sigma _r$ to $\partial\Omega \times \,]-r,r[\,$. 
Then for each coordinate patch $\kappa '_i$,
we can use the mapping $\kappa
_i\colon x'+tn(x')\mapsto (\kappa '_i(x'),t)$ 
as the diffeomorphism in dimension $n$; $\kappa _i$ goes from $U_i$ to
$V_i$, where
$$
U_i=\{x'+tn(x')\mid x'\in U'_i, |t|<r\}, \quad V_i=V'_i\times
\,]-r,r[\,.\tag4.13
$$ 
The advantage is that the normal $\nu (x')$ at $x
'\in\partial\Omega $ is carried over to the normal $(0,1)$ at $(\kappa
'_i(x'),0)$. Moreover, for
points $x\in\Sigma _{r,+}=\Sigma _r\cap \Omega $, $t$ is a
good approximation to the distance function
$d(x)=\operatorname{dist}(x,\partial\Omega )$; their difference goes to
0 for $t\to 0$. 

We can supply these charts with a chart
consisting of the identity mapping on an open set $U_0$ containing $\Omega
\setminus \overline \Sigma _{r,+}$, with $\overline U_0\subset \Omega
$, to get a full cover of $\comega$.
\endexample

Together with the cover by local coordinate charts there exists an associated  partition of unity $\varphi
_0,\dots,\varphi _{I_0}$ such that each $\varphi _i$ is in $ C_0^\infty (U_i)$
taking values in $[0,1]$, and $\sum_{0\le i\le i_0}\varphi _i(x)=1$ for $x\in\comega$. 
It will be
convenient in the following to have the more refined concept of a  partition of unity {\it
 subordinate} to a system of local coordinates, where any two functions are supported
 in one of the $U_i$'s. This fact was originally used in Seeley
 \cite{S69}, proofs are given (in more complicated cases) in
 \cite{G96}, Appendix, and \cite{G09}, Ch.\ 8. For the convenience of the reader we provide a
 proof here.

\proclaim{Lemma 4.4} There exists a system of coordinate charts
$\kappa _i\colon U_i\to V_i$, $i=0,\dots, I_1$, and a subordinate
partition of unity $\varrho _j$, $j=1,\dots,J_0$ (with values in
$[0,1]$ and sum $1$ on $\comega$), such that for each
pair $k,l\le J_0$ there is an $i=i(k,l)\le I_1$ such that
$\operatorname{supp}\varrho _k\cup \operatorname{supp}\varrho
_l\subset U_i$.  
\endproclaim

\demo{Proof} We start out with an arbitrary cover by coordinate charts $\kappa
_i\colon U_i\to V_i$, $i=0,\dots,I_0$. By the
compactness of $\comega$, there is a $\delta >0$ such that any subset
of $\comega$ with diameter $\le\delta $ is contained in one of the
$U_i$'s.
Cover $\comega$ with a finite system of open balls $B_j$ with radius
$\le \delta /4$, $j=1,\dots,J_0$. When $B_{j_1}$ and $B_{j_2}$ are two such
balls, we have two possibilities: 

1) If $B_{j_1}\cap B_{j_2}\ne \emptyset$, it has diameter $\le\delta $, hence
lies in a set $U_i$, take the first such $i$. We shall adjoin the set
$U'=B_{j_1}\cup  B_{j_2}$ to our system, using the mapping $\kappa
_i$ to define a coordinate mapping $\kappa '$ from $U'$ to $V'=\kappa _i(B_{j_1}\cup
B_{j_2})$.  

2) If $B_{j_1}\cap B_{j_2}= \emptyset$, the balls lie in two possibly different sets $U_{i_1}$
and $U_{i_2}$ (take the first $i_1$ and first $i_2$ that occur); then
we shall adjoin
the coordinate neighborhood $U'=B_{j_1}\cup B_{j_2}$ to the given
system using as coordinate transformation the mapping $\kappa _{i_1}$ on $B_{j_2}$  and $\kappa
_{i_2}$ on $B_{j_2}$. Here we may have to make a translation $\tau $ of the
image $\kappa _{i_2}(B_{j_2})$ to make it disjoint from $\kappa
_{i_i}(B_{j_1})$. In this way we get a coordinate chart $\kappa '$ from
$U'$ to $V'=\kappa _{i_1}(B_{j_1})\cup \tau \kappa _{i_2}(B_{j_2})$.

We do this for all pairs $j_1,j_2$ and enumerate the resulting
coordinate charts $\kappa ' \colon U'\to V'$ by numbers $i=I_0+1,\dots,
I_1$;
then we get an extended cover of $\comega$ by coordinate charts
$\kappa
_i\colon U_i\to V_i$, $i=0,\dots,I_1$.

Finally, let $\varrho _j$, $j=1,\dots,J_0$, be a partition of unity
associated with the cover $B_j$, $j=1,\dots, J_0$ (i.e.\ with $\varrho _j\in
C_0^\infty (B_j)$ for each $j$), then any two functions $\varrho _k,
\varrho _l$ have their support in one of the open sets in the extended
cover. 
\qed

\enddemo

We now consider a classical $\psi $do $P$ on $\R^n$ of order $2a$ with
even symbol, elliptic avoiding a ray. It has the $a$-transmission
property with respect to $\Omega $, and an application of Theorem 3.4
in local coordinates shows that the factorization index is $a$. Then
by the general theory of \cite{G15}, the Dirichlet problem (A.10)
satisfies:
When $u\in \dot H^\sigma (\comega)$
 (with $\sigma >a-\frac12$) solves (A.10) for some
 $f\in \ol H^{s-2a}(\Omega )$ with
 $s>a-\frac12$, then $u\in H^{a(s)}(\comega)$;
 moreover, $r^+P$ is Fredholm from  $ H^{a(s)}(\comega)$ to $\ol
 H^{s-2a}(\Omega )$. 
Our principal integration-by-parts theorem is:

\proclaim{Theorem 4.5} Let $P$ be a classical $\psi $do on
${\Bbb R}^n$ of order $2a$ ($0<a<1$), elliptic avoiding a ray, and with even symbol. 
For $u,u'\in H^{a(s)}(\comega)$ with $s\ge a+1$ there
holds, for $j=1,\dots,n$:
$$
\aligned
\int_{\Omega  } Pu\,\partial_j\bar
u'\,dx&+\int_{\Omega }\partial_ju\,\overline{  P^*u'}\,dx\\
&=\Gamma
(a+1)^2\int_{\partial\Omega }s_0\nu _j\gamma _0(d^{-a}u)\,\gamma
_0(d^{-a}\bar u')\, d\sigma+\int_{\Omega } P^{(j)}u\,\bar u'\,dx
 ,\endaligned\tag4.14
$$
where $s_0(x)$ is the value of the principal symbol of $P$ at $(x,\nu
(x))$ for $x\in\partial\Omega $, and $P^{(j)}=P\partial_j-\partial_jP$.

The  term with $P^{(j)}$ vanishes if $P$ is independent of $x_j$ (in particular, when $P$ is 
translation-invariant).

The formula extends to the case $s>a+\frac12$, with the integrals over
$\Omega $ replaced by dualities:
$$
\multline
\ang{r^+ Pu,\partial_j
u'}_{\ol
H^{\frac12-a+\varepsilon }, \dot
H^{a-\frac12-\varepsilon}}+\ang{\partial_ju,{  P^*u'}}_{\dot H^{a-\frac12+\varepsilon},\ol
H^{\frac12-a-\varepsilon }}
\\=\Gamma (a+1)^2\int_{\partial\Omega }\nu _js_0\gamma _0(x_n^{-a}u)\,\gamma _0(x_n^{-a}\bar u')\,
d\sigma +\ang{r^+ P^{(j)}u,
u'}_{\ol
H^{\frac12-a+\varepsilon }, \dot H^{a-\frac12-\varepsilon }};
\endmultline\tag4.15
$$
the last term is a scalar product $(P^{(j)}u,u')_{L_2(\Omega )}$ when
$a\le\frac12$.

\endproclaim

\demo{Proof} For a transparent notation, we formulate the proof in the
case $s\ge a+1$; the extensions to dualities for lower $s$ follow
easily (as in Theorem 3.4).

 Starting with a choice of coordinate charts as in Remark 4.3,
we use Lemma 4.4 to extend it to a covering of $\comega $ with a system of coordinate
patches $\kappa _i\colon U_i\to V_i\subset \R^n$, $i=0,\dots,I_1$,
such that there is a subordinate partition of unity $\varrho _j$, $j=1,\dots , J_0$,
where for any pair of indices $k,l\le J_0$ there is a $U_i$, $i=i(k,l)$,
such that $\varrho _k$ and $\varrho _l$ have support in $U_i$. We can
moreover choose real functions $\psi _k,\psi _l\in C_0^\infty (U_i)$ such that
$\psi _k\varrho _k=\varrho _k$, $\psi _l\varrho _l=\varrho _l$ (i.e.,
they are 1 on the respective supports). Then 
$$
\aligned
\int_{\Omega }( Pu\,\partial_j\bar
u'&+\partial_ju\,\overline{  P^*u'})\,dx=
\sum_{k,l\le J_0}\int_{\Omega }( P\varrho _ku\,\partial_j\bar \varrho _l
u'+\partial_j\varrho _ku\,\overline{  P^*\varrho _lu'})\,dx\\
&=\sum_{k,l\le J_0}\int_{\Omega }( P\psi _k\varrho _ku\,\partial_j \psi _l\varrho _l
\bar u'+\partial_j\psi _k\varrho _ku\,\overline{  P^*\psi _l\varrho
_lu'})\,dx\\
&=\sum_{k,l\le J_0}\int_{\Omega }( P_{kl}u_k\,\partial_j 
\bar u'_l+\partial_ju_k\,\overline{  P_{kl}^*u'_l})\,dx,
\endaligned\tag4.16
$$
where 
$$
P_{kl}=\psi _lP\psi _k,\quad P_{kl}^*=\psi _kP^*\psi _l,\quad
u_k=\varrho _ku,\quad u'_l=\varrho _lu'.\tag4.17
$$
For each pair $(k,l)$ we treat the term by use of the coordinate map for $U_i$,
$i=i(k,l)$.
Denote by $\widetilde P_{kl}$ the operator on $V_i\subset {\Bbb R}^n$ that $P_{kl}$
carries over to; it has compact kernel support in $V_i\times V_i$. 
In detail, $\widetilde P_{kl}=\widetilde \psi _l^{(i)}\widetilde
P^{(i)}\widetilde \psi _k^{(i)}$, cf.\ (4.12).
The
parity  property of the symbol, hence the $a$-transmission property, is preserved
under the coordinate transformation. By Theorem 2.8 applied to $\widetilde P^{(i)}$, $\widetilde
P_{kl}$ has a decomposition into a product of $\pm$-factors and a lower-order term:
$$
\widetilde P_{kl}=\widetilde P_{kl}^-\widetilde P_{kl}^++\widetilde
{ S}_{kl},\text{ in detail }\widetilde P_{kl}^-=\widetilde \psi _l^{(i)}{\widetilde
P^{(i)-}},\widetilde P_{kl}^+={\widetilde
P^{(i)+}}\widetilde \psi _k^{(i)},
\tag4.18
$$
where $\widetilde P_{kl}^\pm$ preserve support in $\crnpm$,
respectively, and $\widetilde {S}_{kl}$ is of order $2a-1$ with
a structure like $S$ in Theorem 3.4, with compact kernel
support in $V_i\times V_i$. We can moreover assume that $\widetilde
P_{kl}^\pm$  have compact kernel supports in $V_i\times V_i$ since
 multiplication by a smooth cutoff function that is 1 on the supports of $\widetilde \psi _k^{(i)}$,
$\widetilde \psi _l^{(i)}$, changes the operator by a smoothing term.

Now all this is carried back to $U_i$ by the coordinate
transformation; $\widetilde P_{kl}^\pm$ are carried over to operators
$P_{kl}^\pm$, and $\widetilde {\Cal S}_{kl}$ is carried over to $\Cal S_{kl}$. The property that $\widetilde P_{kl}^\pm$ preserve
supports in $\crnpm$, respectively, carries over to the property that 
$P_{kl}^\pm$ preserve support in $\comega$ resp.\ $\complement\Omega
$.
Then we have the adjoint mapping properties (where $r^+$ and $e^+$ are defined
relative to $\comega\subset {\Bbb R}^n$):
$$\aligned
r^+P^-_{kl}e^+&\colon\ol H_{p}^s(\Omega )\to \ol H_{p}^{s-a}(\Omega )\text{ and
 }{P_{kl}^-}^*\colon \dot H_{p'}^{a-s}(\comega )\to \dot
H_{p'}^{-s}(\comega )\text{ are adjoints},\\
r^+{P_{kl}^+}^*e^+&\colon\ol H_{p}^s(\Omega )\to\ol H_{p}^{s-a}(\Omega )\text{ and
 }P_{kl}^+\colon \dot H_{p'}^{a-s}(\comega )\to \dot
 H_{p'}^{-s}(\comega )\text{
 are adjoints}.
\endaligned \tag4.19
$$

With this preparation, we can calculate as follows: Denote
 $r^+P_{kl}^+u_k=w$, $r^+{P_{kl}^-}^*u_l'=w'$.
Then
$$
I=\int_{\Omega \cap U_i}( P_{kl}^-P_{kl}^+u_k\,\partial_ju_l'
+\partial_ju_k\,\overline{ {P_{kl}^+}^*{P_{kl}^-}^*u_l'})\, dx=
\int_{\Omega \cap U_i}( P_{kl}^-e^+w\,\partial_ju_l'
+\partial_ju_k\,\overline{ {P_{kl}^+}^*e^+w'})\, dx.
$$
We apply Theorem 4.2 to the first term, and a conjugated variant
to the second term, obtaining
$$\aligned
I&=2\int_{\partial\Omega \cap U_i}\nu _j\,\gamma _0w\,\gamma _0\bar
w'\, d\sigma +\int_{\Omega \cap U_i}(w\,\partial_j\bar w'+\partial_jw\,\bar w'
+w\,\overline{[{P_{kl}^{-}}^*,\partial_j]u_l'}
+[P_{kl}^+,\partial_j]u_k\,\bar w')\, dx\\
&=\int_{\partial\Omega \cap U_i}\nu _j\,\gamma _0w\,\gamma _0\bar
w'\, d\sigma +\int_{\Omega \cap U_i} [P_{kl}^-e^+r^+P_{kl}^+,\partial_j]u_k\,\bar u_l'\, dx.
\endaligned$$
For the second line it was used that $\int_{\Omega \cap
U_i}(w\,\partial_j\bar w'+\partial_jw\,\bar w')\, dx'$ gives another
copy of $\int_{\partial\Omega \cap U_i}\nu _j\,\gamma _0w\,\gamma _0\bar
w'\, d\sigma $ with a minus sign, and the two terms with commutators
were reduced to a single term as in the proof of Theorem 3.4.
  
For the term with $S_{kl}$ we proceed as in Theorem 3.4, concluding
that it gives  no
boundary contribution, only a commutator term that can be added to the one
with $P_{kl}^-e^+r^+P_{kl}^+$. 

This leads to the formula
$$
\int_{\Omega }( P_{kl}u_k\,\partial_j 
\bar u'_l+\partial_ju_k\,\overline{  P_{kl}^*u'_l})\,dx=\int_{\partial\Omega \cap U_i}\nu _j\,\gamma _0w\,\gamma _0\bar
w'\, d\sigma +\int_{\Omega \cap U_i} [P_{kl},\partial_j]u_k\,\bar
u_l'\, dx.\tag 4.20
$$

The boundary contributions from $P_{kl}^\pm$ are found from the values of the
functions in the localized situation. Here $\gamma _0(P_{kl}^+u_k)$
comes from $$
\widetilde P^{(i)+}\widetilde\psi _{k}^{(i)}\widetilde{\varrho
_ku}^{(i)}|_{x_n=0}=\widetilde P^{(i)+}\widetilde{\varrho
_ku}^{(i)}|_{x_n=0}=\lim_{x_n\to 0+}\Gamma (a+1)x_n^{-a}\widetilde{\varrho
_ku}^{(i)},
$$
by calculations as in (3.31); recall that $\psi _k\varrho _k=\varrho _k$. This carries over to $\partial\Omega $
as \linebreak$\Gamma (a+1)\lim_{d\to 0}(d^{-a}\varrho _ku)$, since $\widetilde
d^{(i)}/x_n\to 1$ for $x_n\to 0$. Similarly, cf.\ (3.32), $\gamma
_0({P_{kl}^-}^*u'_l)= \Gamma (a+1)s_0\gamma _0(d^{-a}\varrho _lu')$.
We conclude that
$$
 \int_{\partial\Omega \cap U_i}\nu _j\gamma _0(P_{kl}^+u_k)\gamma _0(\overline{{P_{kl}^{-}}^*u_l'}
)\, d\sigma =\Gamma (a+1)^2 \int_{\partial\Omega \cap U_i}\nu
_js_0\gamma _0(d^{-a}\varrho _ku)\gamma _0(d^{-a}\varrho _l\bar u')\, d\sigma .
$$

We have then obtained: 
$$
\aligned
&\int_{\Omega \cap U_i}( P_{kl}u_k\,\partial_j\bar u_l'
+\partial_ju_k\,\overline{ P_{kl}^*u_l'})\, dx\\
&=\Gamma (a+1)^2 \int_{\partial\Omega \cap U_i}\nu
_js_0\gamma _0(d^{-a}\varrho _ku)\gamma _0(d^{-a}\varrho _l\bar u')\, d\sigma + \int_{\Omega \cap U_i} [P_{kl},\partial_j]u_k\,\bar u_l'\, dx,
\endaligned
$$
for each pair $(k,l)$, and when we sum over $k$ and $l$,
using that $\sum _k\varrho _k=\sum_l\varrho _l=1$ on $\comega$, we
find (4.14). 

The extension to dualities in (4.15), when $s>a+\frac12$, follows when one
formulates the detailed study of $P_{kl}$ in terms of dualities as in
Theorem 3.4.
\qed
\enddemo

The validity extends to suitable H\"older spaces. To get a very efficient statement, we can apply the
general result of \cite{G14} Th.\ 4.2, Ex.\ 4.3, for H\"older-Zygmund
spaces, showing that $r^+P$ defines a Fredholm operator for $s>a-1$:
$$
r^+P\colon C_*^{a(s)}(\comega)\to \ol C_*^{s-2a}(\Omega). \tag4.21
$$
There is also a regularity result stating that when $u\in \dot
 C_*^{t}(\comega )$ for some $t>a-1$ (in particular when $u\in e^+L_\infty (\Omega)$), then $r^+Pu\in
\ol C_*^{s-2a}(\Omega)$ implies $u\in  C_*^{a(s)}(\comega)$.
We recall that $\ol C^s_*(\Omega)$ equals the H\"older space
$C^s(\comega)$ when $s>0$, $s\notin{\Bbb N}$; cf.\ also (A.3).
Here  $C_*^{a(s)}=
\Lambda _+^{(-a)}e^+\ol C_*^{s-a}(\Omega)$, where $\Lambda _+^{(t)}$ is an
order-reducing operator on ${\Bbb R}^n$ 
preserving support in
$\comega$, as recalled in (A.7) and used in the proof of Theorem 4.1. These operators
apply also to $C^s_*$-spaces by \cite{G14}. 

To assure that $r^+Pu$ is bounded and  $\partial_ju$ is integrable 
on $\Omega $, we take $s=1+a+\varepsilon $ 
with $\varepsilon >0$. 
Then $r^+Pu\in \ol C^{1-a+\varepsilon }(\Omega)$, and 
(when $1+a+\varepsilon \notin{\Bbb N}$)
$$
u\in C_*^{a(1+a+\varepsilon )}(\comega)\subset e^+d^a\ol C^{1+\varepsilon }(\Omega),\tag4.22
$$
with $\partial_ju\in e^+d^{a-1}\ol C^{1+\varepsilon
}(\Omega)+e^+d^a\ol C^\varepsilon (\Omega)
\subset L_1(\Omega )$. Since the various spaces are invariant under
$C^\infty $-coordinate changes, the proof of Theorem 4.5 carries
through
for such functions.

We have hereby obtained:

\proclaim{Corollary 4.6} 
 Formula {\rm (4.14)}
holds also when $u,u'\in C_*^{a(1+a+\varepsilon )}(\comega)$,
 some $\varepsilon >0$.

This is assured when $u,u'\in e^+L_\infty (\Omega)$
 and $r^+Pu, r^+P^*u'\in \ol C^{1-a+\varepsilon '}(\comega)$
($\varepsilon '=\varepsilon $ when $1+a+\varepsilon \notin{\Bbb N}$,
 $\varepsilon '>\varepsilon $ when $1+a+\varepsilon \in{\Bbb N}$).
\endproclaim

The assumption on $r^+Pu$ in the corollary is a little more general than the assumption in \cite{RS14b,
RSV15} which take $r^+Pu\in C^{0,1}(\comega)$.  On the other hand, these authors work under a weaker smoothness
hypothesis on $\Omega $ (namely that it is $C^{1,1}$).

The assumptions in Theorm 4.5 are a considerable generalization.

The advantage of referring to $H^{a(s)}(\comega)$ and $C^{a(s)}_*(\comega)$
is that these scales of spaces do not depend on a choice of $P$, but are
the appropriate solution spaces for the Dirichlet problem for  {\it all}
classical elliptic $\psi $do's  $P$ of order $2a$ and type $a$ with
factorization index $a$.  

The results apply for example to $(-\Delta )^a$ 
 and to $a$'th powers $A^a$ of
second-order strongly elliptic differential operators $A$ with $C^\infty $-coefficients. 
Seeley \cite{S67} showed that $A^a$ is a classical $\psi $do of order $2a$,
 with a symbol constructed via the resolvent; it is even. $A^a$ is
 again strongly elliptic, since the principal symbol is $(a_0(x,\xi
 ))^a$, taking values in $\{\operatorname{Re}z>0\}$ for $|\xi |\ge1$.
For $(-\Delta
 )^a$ one can more directly remark that the symbol may be written
$|\xi |^{2a}=|\xi |^{2a}\eta (\xi )+  |\xi |^{2a}(1-\eta (\xi ))
$ with an excision function $\eta $ (cf.\ (2.11)), and proceed in a
 similar
 way as in Corollary 3.5.

As a consequence of the above results, we can moreover show
an integration-by-parts formula where $\partial_j$ is replaced by a
radial derivative $x\cdot\nabla=\sum_{j=1}^nx_j\partial_j$.

\proclaim{Theorem 4.7} Let $P$ be a classical elliptic $\psi $do on
${\Bbb R}^n$ of order $2a$ ($0<a<1$) with even symbol. 
Then for $u,u'$ as in Theorem {\rm 4.5} or Corollary {\rm 4.6} there
holds:
$$
\align
\int_{\Omega  }( Pu\,(x\cdot \nabla \bar u')+(x\cdot
\nabla u)\, \overline {P^*u'})\,dx
&=\Gamma
(a+1)^2\int_{\partial\Omega }(x\cdot\nu)s_0 \gamma _0(d^{-a}u)\,\gamma
_0(d^{-a}\bar u')\,
d\sigma\\
&\quad -n\int_{\Omega } Pu\,\bar u'\,dx+\int_\Omega
[P,x\cdot\nabla]u\,\bar u'\,dx;\tag4.23
\endalign
$$
here
$$
[P,x\cdot\nabla]=P_1-P_2,\quad P_1=\Op(\xi \cdot\nabla_\xi p(x,\xi
)),\quad P_2=\Op(x \cdot\nabla_x p(x,\xi )).\tag4.24
$$

When
$u\in H^{a(s )}(\comega)$ with $a+\frac12<s< a+1$, some integrals  are replaced by dualities:
$$
\multline
\ang{r^+ Pu,x\cdot \nabla u'}_{\ol
H^{\frac12-a+\varepsilon }, \dot
H^{a-\frac12-\varepsilon}}+\ang{x\cdot \nabla u,{  P^*u'}}_{\dot H^{a-\frac12+\varepsilon},\ol
H^{\frac12-a-\varepsilon }}
\\
=\Gamma (a+1)^2((x\cdot\nu )s_0\gamma _0(d^{-a}u),\gamma _0(d^{-a}u'))_{L_2(\R^n)}\\ -n\ang{r^+ Pu,
u'}_{\ol
H^{\frac12-a+\varepsilon }, \dot H^{a-\frac12-\varepsilon }}
+
\ang{r^+ [P,x\cdot\nabla]u,
u'}_{\ol
H^{\frac12-a+\varepsilon }, \dot H^{a-\frac12-\varepsilon }}.
\endmultline\tag4.25
$$
\endproclaim

\demo{Proof} The calculation goes as follows:
$$
\aligned
&\int_{\Omega  } Pu\,(x\cdot \nabla \bar u')\,dx+
\int_{\Omega  }(x\cdot \nabla
u)\, \overline{P^*u'}\,dx\\
&=\sum_{j=1}^n\int_{\Omega
}(x_j Pu\,\partial_j\bar
u'+\partial_j(x_ju)\, \overline{P^*u'}-u\, \overline{P^* u'}
)\,dx\\
&=\sum_{j=1}^n\int_{\Omega  }(
 P(x_ju)\,\partial_j\bar u'+ [x_j,P]u\,\partial_j\bar
u'+\partial_j(x_ju)\, \overline{P^*u'}- Pu\,\bar
u')\,dx\\
&=\Gamma
(a+1)^2\int_{\partial\Omega }( x\cdot\nu )\,s_0\,\gamma _0(d^{-a}u)\gamma
_0(d^{-a}\bar u')\,
d\sigma -n \int_{\Omega }Pu\,\bar u')\,dx
\\
&\quad+\int_{\Omega }\sum_j [P,\partial_j]x_ju\,\bar u'
 +\int_\Omega \sum_j [x_j,P]u\,\partial_j\bar u'\,dx 
\endaligned\tag4.26
$$
For the second equality we have applied Theorem 4.1 to $u\,
\overline{P^* u'}$, and for the third equality
we have applied Theorem 4.5 to the terms $ P(x_ju)\,\partial_j\bar
u'$ and $\partial_j(x_ju)\, \overline{P^*u'}$.

For the last term, we observe that  $[x_j,P]$ equals $\Op(i\partial_{\xi _j}p(\xi ))$, which is a classical $\psi $do of order
$2a-1$, again with even symbol (having the $a$-transmission
property at $\partial\Omega $), so $[x_j,P]u\in H^{\frac32-a+\varepsilon }(\Omega
)$ resp.\ $C^{2-a+\varepsilon
}(\comega)$, and $\partial_j[x_j,P]u\in H^{\frac12-a+\varepsilon
}(\Omega )$ resp.\ $ C^{1-a+\varepsilon
}(\comega)$, under the hypotheses in Theorem 4.5 
resp.\ Corollary  4.6 (by  \cite{G15}, Th.\ 4.2, resp.\ \cite{G14},
Th.\ 3.2(1)). Then
$$
\int_\Omega  [x_j,P]u\,\partial_j\bar u'\,dx+\int_\Omega
 \partial_j[x_j,P]u\,\bar u'\,dx=\int_{\partial \Omega }\nu _j\gamma
_0([x_j,P]u)\gamma _0\bar u'\,dx =0,\tag4.27
$$
since $\gamma _0u'=0$, so $\int_\Omega  [x_j,P]u\,\partial_j \bar u'\,dx=-\int_\Omega
 \partial_j[x_j,P]u\,\bar u'\,dx$. 

Moreover,
$$
[P,\partial_j]x_ju-\partial_j[x_j,P]u=P\partial_jx_ju-\partial_jx_jPu=[P,x_j\partial_j]u,
$$
so the two commutator integrals with $\partial_j$ and $x_j$ together give $\int_{\Omega
} [P,x\cdot\nabla]u\,\bar u'\,dx$. This shows (4.23).

Considering the symbols, since $[P,\partial_j]$ has symbol $-\partial_{x_j}p$
and $[x_j,P]$ has symbol $i\partial_{\xi _j}p$ (by the formula for the
Leibniz product,
cf.\ (2.34)),
$$
\multline
\operatorname{symbol}([P,\partial_j]x_j-\partial_j[x_j,P])=-\partial_{x_j}p\,\#\,
x_j-i\xi _j\,\# \,i\partial_{\xi _j}p\\
=-x_j\partial_{x_j}p-(-i)\partial_{\xi _j}\partial_{x_j}p+\xi _j\partial_{\xi _j}p+(-i)\partial_{x_j}\partial_{\xi _j}p=-x_j\partial_{x_j}p+\xi _j\partial_{\xi _j}p,
\endmultline
$$
so $ [P,x\cdot\nabla]$ has symbol $\xi \cdot\nabla_\xi
p(x,\xi )-x\cdot\nabla_x p(x,\xi )$; this shows (4.24).

The result extends to spaces with lower $s$ as in Theorem 4.5, in the form  (4.25).
\qed 
\enddemo

Observe some special cases:

\proclaim{Corollary 4.8} In the situation of Theorem {\rm 4.7},
if $P$ is $x$-independent then $P_2$ vanishes. If, in addition, the
symbol $p$ of $P$ is homogeneous of degree $2a$ (i.e., equals it
principal part), then $P_1=2aP$, and formula {\rm (4.23)} takes the form
$$
\align
\int_{\Omega  }( Pu\,(x\cdot \nabla \bar u')+(x\cdot
\nabla u)\, \overline {P^*u'})\,dx
&=\Gamma
(a+1)^2\int_{\partial\Omega }(x\cdot\nu)s_0 \gamma _0(d^{-a}u)\,\gamma
_0(d^{-a}\bar u')\,
d\sigma\\
&\quad +(2a -n)\int_{\Omega } Pu\,\bar u'\,dx.\tag4.28
\endalign
$$
\endproclaim
\demo{Proof} The first statement is obvious. For the second statement,
Euler's formula gives that $\xi \cdot\nabla_\xi p=2ap$, hence
$P_1=2aP$, and the formula follows by insertion. \qed
\enddemo

Formula (4.28) for $P=(-\Delta )^a$ (and real $u=u'$)  was a
principal result of \cite{RS14b}, and was extended to
selfadjoint positive homogeneous $x$-independent  operators $P$ in
\cite{RSV15}, under lower smoothness assumptions that ours.  
 It leads to a Pohozaev-type formula
(generalizing a formula of Pohozaev \cite{P65} for $\Delta $) that
can be used to obtain uniqueness and (non)existence results. We similarly
find from (4.23):

\proclaim{Corollary 4.9} Let $P$ be as in Theorem {\rm 4.7} and selfadjoint, and
let $u$ be a bounded real solution of the problem
$$
 r^+Pu=f(u)\text{ in }\Omega ,\quad \operatorname{supp}u\subset \overline\Omega ,\tag4.29
$$
where $f$ is a real $C^{0,1}$-function. Let $F(t)=\int_0^tf(s)\,ds$. Then
$$
-2n\int_\Omega F(u)\,dx+n\int_\Omega
f(u)\, u\, dx=\Gamma (1+a)^2\int_{\partial\Omega }(x\cdot \nu) \,s_0
\gamma _0(d^{-a}u)^2
\, d\sigma+\int_{\Omega }[P,x\cdot\nabla]u\,u\,dx, \tag4.30
$$
where $[P,x\cdot\nabla]=P_1-P_2$ as in  {\rm (4.24)}.

If $P$ is $x$-independent, the formula becomes
$$
-2n\int_\Omega F(u)\,dx+n\int_\Omega
f(u)\, u\, dx=\Gamma (1+a)^2\int_{\partial\Omega }(x\cdot \nu) \,s_0
\gamma _0(d^{-a}u)^2
\, d\sigma+\int_{\Omega }P_1u\,u\,dx. \tag4.31
$$
Here
if the symbol of $P$ 
moreover is homogeneous,  the formula reduces to
$$
-2n\int_\Omega F(u)\,dx+(n-2a)\int_\Omega
f(u)\, u\, dx=\Gamma (1+a)^2\int_{\partial\Omega }(x\cdot \nu) \,s_0
\gamma _0(d^{-a}u)^2
\, d\sigma. \tag4.32
$$ 
\endproclaim

\demo{Proof} Since $P=P^*$, the left-hand side of (4.23) reduces for
real $u=u'$ to \linebreak$2\int_{\Omega }Pu\, (x\cdot\nabla u)\, dx$. Since $u$ is bounded, so is $f(u)$; then $u\in \dot
C^a(\comega)$ by the regularity theory. Then since $F(0)=0$, $F(u)\in
\dot C^a(\comega)$. We have that
$$
\aligned
(x\cdot\nabla)F(u)&=\sum_{j=1}^nx_j\partial_jF(u)=\sum_{j=1}^nx_jF'(u)\partial_ju=f(u)(x\cdot\nabla u),\\
(x\cdot\nabla)F(u)
&=\sum_{j=1}^n\partial_j(x_jF(u))-nF(u).
\endaligned
$$
Then since the integral over $\Omega $ of $\partial_j(x_jF(u))$ is zero,
$$
\int_{\Omega }(x\cdot\nabla u)f(u)\,dx=\int_{\Omega }(x\cdot\nabla )F(u)\,dx=-n\int_{\Omega }F(u)\,dx.
$$
Insertion of this and the formula $f(u)=r^+Pu$ in (4.23) leads to
(4.30). 

The last statements follow as in Corollary 4.8.
\qed

\enddemo

Formula (4.32) is the formula shown in \cite{RSV15}.

The new result will for example apply to fractional powers of  magnetic
Schr\"odinger operators. To draw conclusions on solvability of
nonlinear equations, one will have to investigate sign properties of
the involved integrals.

Let us just end here by illustrating the use in some very simple examples in the
$x$-independent case.

\example{Example 4.10} The fractional Helmholtz (or Schr\"odinger) 
operator 
$P=(-\Delta +m^2)^a$, $0<a<1$ and $m>0$, has
the symbol $p(\xi )=(|\xi |^2+m^2)^a$ of order $2a$. It is not
homogeneous, but has the (classical)  expansion in homogeneous terms
$$
p(\xi )\sim |\xi |^{2a}+am^2|\xi |^{2a-2}+\tfrac12 a(a-1)m^4|\xi |^{2a-4}+\dots,
$$
and it is even. In this case  $$
\xi
\cdot\nabla p(\xi )=2a|\xi |^2(|\xi |^2+m^2)^{a-1}>0\text{ for }\xi \ne
0,$$
and $P_1=\operatorname{Op}(\xi
\cdot\nabla p(\xi ))=2a(-\Delta )(-\Delta +m^2)^{a-1}$ is  positive
on the functions $v$ in the Dirichlet domain $H^{a(2a)}(\comega)$ of
$P$: Here $v\in \dot H^a(\comega)\subset H^a(\R^n)$ and $r^+P_1v\in
L_2(\Omega )$, $P_1v\in H^{-a}(\R^n)$, so
$$
 \int_{\Omega
}r^+P_1v\,\bar v\,dx=\ang{P_1v,v}_{H^{-a}(\R^n),H^a(\R^n)}=\tfrac{2a}{(2\pi )^{n}}\int_{\R^n}|\xi
|^2(|\xi |^2+m^2)^{a-1}|\hat v(\xi )|^2\, d\xi >0
$$
unless $v\equiv 0$.
Let us see
what this gives for an eigenvalue problem
$$
r^+Pu=\lambda u\text{ in }\Omega ,\quad \operatorname{supp}u\subset \overline\Omega ,\tag 4.33
$$
   for some $\lambda \in {\Bbb R}$ and bounded  real $u$. With $f(u)=\lambda u$, $F(u)=\frac12
   \lambda u^2$, so the first two integrals in (4.31) cancel out, giving
$$
0=\Gamma (1+a)^2\int_{\partial\Omega }(x\cdot \nu) \,s_0
\gamma _0(d^{-a}u)^2
\, d\sigma+ \int_{\Omega }P_1u\,u\,dx.
$$
By the positivity of $P_1$, this allows the conclusion
$$
\gamma _0(d^{-a}u)=0\implies \int_{\Omega
}P_1u\,u\,dx=0\implies u\equiv 0.
$$
This shows a kind of {\it unique continuation principle} for solutions of
the eigenvalue equation: When $u$ is in the Dirichlet domain and in
addition $\gamma _0(d^{-a}u)=0$, then $u\equiv 0$.

\endexample

\example{Example 4.11} For the operator in Example 4.10,
$$
\aligned
P_1&=2a(-\Delta )(-\Delta +m^2)^{a-1}=2a(-\Delta +m^2)^a-2am^2(-\Delta
+m^2)^{a-1}=2aP-P_3,\\
&\text{ with }P_3=2am^2(-\Delta
+m^2)^{a-1},
\endaligned
$$
here $P_3$ is a positive operator. Thus for bounded real solutions of
(4.29),
equation (4.31) can be written in the form
$$
-2n\int_\Omega F(u)\,dx+(n-2a)\int_\Omega
f(u)\, u\, dx+\int_{\Omega }P_3u\,u\,dx=\Gamma (1+a)^2\int_{\partial\Omega }(x\cdot \nu) \,s_0
\gamma _0(d^{-a}u)^2
\, d\sigma. \tag4.34
$$
Consider the case $f(u)=u|u|^{r-1}=\operatorname{sign}u\,|u|^r$ with an $r>1$. Here since
$F(u)=\frac1{r+1}|u|^{r+1}$, (4.34) takes the form
$$
\tfrac{-2n+(n-2a)(r+1)}{r+1}\int_\Omega |u|^{r+1}\,dx+\int_{\Omega }P_3u\,u\,dx
=\Gamma (1+a)^2\int_{\partial\Omega }(x\cdot \nu) \,s_0
\gamma _0(d^{-a}u)^2
\, d\sigma.
\tag4.35
$$
 Consider a starshaped domain $\Omega $
 ($n\ge 2$); we can assume that $0$ is a center. Then $x\cdot\nu \le
 0$ on $\partial\Omega $ (recall that our $\nu $ is the interior
 normal).
Note that
$$
[-2n+(n-2a)(r+1)=]\;(n-2a)r-(n+2a)\gtreqless 0\iff r\gtreqless
 \tfrac{n+2a}{n-2a}.
$$
In the critical and supercritical cases $r\ge \frac{n+2a}{n-2a}$ we thus have that
if $u$ is a bounded solution (hence  is in $\dot C^a(\comega)$), then the left-hand side of
(4.35) is $>0$ unless $u\equiv 0$, and the right-hand side is $\le
0$.

This shows {\it nonexistence of nontrivial  solutions, when  $r\ge
\frac{n+2a}{n-2a}$}.

There is a treatment of existence questions in \cite{RS15a}, which
goes beyond the case of homogeneous integral operator 
kernels, by allowing nonnegative kernels with certain growth estimates
on rays. That approach may possibly also be applicable to this
example.

\endexample

\head{Appendix: Spaces and pseudodifferential operators}\endhead

We here collect the notation and concepts from the theory of
pseudodifferential operators that will be used, including some results from
\cite{G15,G14}. Since the set-up is explained in a much more elaborate
form there, in particular in \cite{G15}, we shall just give a brief
summary here. 

A pseudodifferential operator ($\psi $do) $P$ on ${\Bbb R}^n$ is
defined from a symbol $p(x,\xi )$ on ${\Bbb R}^n\times{\Bbb R}^n$ by 
$$
Pu=p(x,D)u=\operatorname{Op}(p(x,\xi ))u 
=(2\pi )^{-n}\int e^{ix\cdot\xi
}p(x,\xi )\hat u\, d\xi =\Cal F^{-1}_{\xi \to x}(p(x,\xi )\hat u(\xi
));\tag A.1
$$  
here $\Cal F$ is the Fourier transform $(\F u)(\xi )=\hat u(\xi
)=\int_{{\Bbb R}^n}e^{-ix\cdot \xi }u(x)\, dx$. 
We refer to
textbooks such as H\"ormander \cite{H85}, Taylor \cite{T91}, Grubb \cite{G09} for the rules of
calculus. \cite{G09} moreover gives an account of the Boutet de Monvel
calculus of pseudodifferential boundary problems, cf.\ also e.g.\
\cite{G96}. A standard choice is
to take $p$ in the symbol space $S^r_{1,0}({\Bbb R}^n\times{\Bbb R}^n)$, consisting of
$C^\infty $-functions $p(x,\xi )$
such that $\partial_x^\beta \partial_\xi ^\alpha p(x,\xi
)$ is $O(\ang\xi ^{r-|\alpha |})$ for all $\alpha ,\beta $, for some
$r\in{\Bbb R}$; then $p$ and $P$ have order $r$.
Also more general symbol spaces will be used in this paper.
When $P$ is a $\psi $do on ${\Bbb R}^n$, $P_+=r^+Pe^+$
denotes its truncation to $\rnp$, or to $\Omega $, depending on the context.

Let $1<p<\infty $  (with $1/p'=1-1/p$), then the $L_p$-Sobolev spaces
(Bessel-potential spaces) are defined for $s\in{\Bbb R}$ by
$$
\aligned
H^s_p(\R^n)&=\{u\in \SD'({\Bbb R}^n)\mid \F^{-1}(\ang{\xi }^s\hat u)\in
L_p(\R^n)\},\\
\dot H^{s}_p(\comega)&=\{u\in H^{s}_p({\Bbb R}^n)\mid \supp u\subset
\comega \},\\
\ol H^{s}_p(\Omega)&=\{u\in \D'(\Omega )\mid u=r^+U \text{ for some }U\in
H^{s}_p(\R^n)\};
\endaligned \tag A.2
$$
here $\operatorname{supp}u$ denotes the support of $u$. The definition
is also used with $\Omega =\rnp$. In most current texts, $\ol
H^s_p(\Omega )$ is denoted $H^s_p(\Omega )$ without the overline (that
was introduced along with the notation $\dot H$ in \cite{H65,H85}), but we keep it here since it is
practical in indications of dualities, and makes the notation more
clear in formulas where
both types occur. When $p=2$, the mention of
$p$ is usually left out.

We recall that $\ol H_p^s(\Omega )$ and $\dot H_{p'}^{-s}(\comega)$ are dual
spaces with respect to a sesquilinear duality extending the $L_2(\Omega )$-scalar
product, written e.g.\
$$
\ang{f,g}_{\ol H_p^s(\Omega ),\dot H_{p'}^{-s}(\comega)},\text{ or just }\ang{f,g}_{\ol H_p^s,\dot H_{p'}^{-s}}.
$$

There is a wealth of other interesting scales of spaces,
the
Triebel-Lizorkin and Besov  spaces $F^s_{p,q}$ and $B^s_{p,q}$, where
the problems can be studied; see details in \cite{G14}. In the present
work, we shall just use the  H\"older-Zygmund spaces $B^s_{\infty
\infty }$, also denoted $C^s_*$. These are interesting because $C^s_*({\Bbb R}^n)$ equals
the H\"older space $C^s({\Bbb R}^n)$ when $s\in\rp\setminus {\Bbb
N}$. There are similar
statements for derived spaces over $\rnp$ and $\Omega $, and again the
conventions $\ol C$ and $\dot C$ are used for spaces of restricted resp.\ supported functions. For integer values one has, with $C^k_b({\Bbb R}^n)$ denoting the space of functions with bounded
continuous derivatives up to order $k$,
$$
\aligned
C^k_b({\Bbb R}^n)&\subset C^{k-1,1}({\Bbb R}^n)\subset C^k_*({\Bbb R}^n)\subset C^{k-0}({\Bbb R}^n)\text{ when }k\in
{\Bbb N},\\
 \quad C^0_b({\Bbb R}^n)&\subset L_\infty ({\Bbb R}^n)\subset
 C^0_*({\Bbb R}^n),
\endaligned\tag A.3
$$
and similar statements for derived spaces.

We  use the notation  $\bigcup_{\varepsilon >0}
H_p^{s+\varepsilon }(\R^n) =
H_p^{s+ 0}(\R^n)$, $\bigcap_{\varepsilon >0} H_p^{s-\varepsilon }(\R^n)=
H_p^{s- 0}(\R^n)$,
 applied in a similar way for the other scales of
spaces.

A $\psi $do $P$ is called classical (or polyhomogeneous) when the
symbol $p$ has an asymptotic expansion $p(x,\xi )\sim \sum_{j\in{\Bbb
N}_0}p_j(x,\xi )$ with $p_j$ homogeneous in $\xi $ of degree $m-j$ for
all $j$, and $p(x,\xi )- \sum_{j<J}p_j(x,\xi )\in S^{m-J}_{1,0}({\Bbb
R}^n\times \R^n)$ for all $J$. Then $P$ has order $m$. One can even
allow $m$ to be complex (with complex homogeneities, $p_j(x,t\xi
)=t^{m-j}p(x,\xi )$ for $|\xi |\ge 1$, $t\ge 1$);
then $p
$ and its remainders are in 
$ S^{\operatorname{Re} m-J}_{1,0}({\Bbb R}^n\times{\Bbb R}^n)$;
the operator and symbol are still said to be of order $m$. 

Here there is an additional definition, introduced by H\"ormander in
\cite{H65,H85}: $P$ satisfies {\it the $\mu
$-transmission condition} at $\partial\Omega $ (in short: {\it is of type} $\mu $) for some $\mu
\in{\Bbb C}$  when, in local coordinates,
$$
\partial_x^\beta \partial_\xi ^\alpha p_j(x,-\nu )=e^{\pi i(m-2\mu -j-|\alpha | )
}\partial_x^\beta \partial_\xi ^\alpha p_j(x,\nu ),\tag A.4
$$
for all $x\in\partial\Omega $, all $j,\alpha ,\beta $, where 
$\nu $ denotes the interior normal to $\partial\Omega $ at $x$.
The implications of the $\mu $-transmission condition 
were a main subject of \cite{G15}.
  
A special role in the theory is played by the {\it order-reducing
operators}. There is a simple definition of operators $\Xi _\pm^\mu $ on
${\Bbb R}^n$
$$ 
\Xi _\pm^\mu =\operatorname{Op}(([\xi ']\pm i\xi _n)^\mu )\tag A.5 
$$
(or with $[\xi ']$ replaced by $\ang{\xi '}$); they preserve support
in $\crnpm$, respectively. Here the function
$([\xi ']\pm i\xi _n)^\mu $ does not satisfy all the estimates
required for the class $S^{\operatorname{Re}\mu }_{1,0}({\Bbb
R}^n\times{\Bbb R}^n)$, but the operators are useful for many
purposes. There is a more refined choice $\Lambda _\pm^\mu $
\cite{G90,G15}, with
symbols $\lambda _\pm^\mu (\xi )$ that do
satisfy all the estimates for $S^{\operatorname{Re}\mu }_{1,0}({\Bbb
R}^n\times{\Bbb R}^n)$; here $\overline{\lambda _+^\mu }=\lambda _-^{\overline\mu }$.
The symbols have holomorphic extensions in $\xi _n$ to the complex
halfspaces ${\Bbb C}_{\mp}=\{z\in{\Bbb C}\mid
\operatorname{Im}z\lessgtr 0\}$, and hence the operators preserve
support in $\crnpm$, respectively; operators with that property are
called plus- resp.\ minus-operators. There is also a pseudodifferential definition $\Lambda
_\pm^{(\mu )}$ adapted to the situation of a smooth domain $\Omega
$. 

It is elementary to see by the definition of the spaces $H^s_p(\R^n)$
in terms of Fourier transformation, that the operators define homeomorphisms 
for all $s$:
$$
\Xi^\mu _\pm\colon H^s_p(\R^n) \simto H^{s-\operatorname{Re}\mu
}_p(\R^n), \quad  
\Lambda ^\mu _\pm\colon H^s_p(\R^n) \simto H^{s-\operatorname{Re}\mu
}_p(\R^n)\tag A.6
$$
(and so does of course $\Xi ^\mu =\Op(\ang \xi ^\mu )$). The special
interest is that the plus/minus operators also 
 define
homeomorphisms related to $\crnp$ and $\comega$: 
$$
\aligned
\Xi ^{\mu }_+, \Lambda ^\mu _+\colon \dot H^s_p(\crnp )\simto
\dot H^{s-\operatorname{Re}\mu }_p(\crnp),&\quad
\Xi ^{\mu }_{-,+}, \Lambda ^\mu _{-,+}\colon \ol H^s_p(\rnp )\simto
\ol H^{s-\operatorname{Re}\mu }_p(\rnp );
\\
\Lambda^{(\mu )}_+\colon \dot H^s_p(\comega )\simto
\dot H^{s-\operatorname{Re}\mu }_p(\comega ),&\quad
\Lambda ^{(\mu )}_{-,+}\colon \ol H^s_p(\Omega )\simto
\ol H^{s-\operatorname{Re}\mu }_p(\Omega );
\endaligned \tag A.7
$$
for all $s\in{\Bbb R}$;
 here $\Xi ^\mu _{-,+}$, $\Lambda ^\mu _{-,+}$ resp.\  $\Lambda ^{(\mu )}_{-,+}$ is short for
$r^+\Xi ^\mu _-e^+$,  $r^+\Lambda  ^\mu _-e^+$  
resp.\  $r^+\Lambda ^{(\mu
 )}_{-}e^+$, suitably extended to large negative $s$ (cf.\ Rem.\ 1.1
 and Th.\ 1.3 in \cite{G15}). 

One has moreover, that
 the operators $\Xi ^\mu _{+}$ and $r^+\Xi ^{\overline\mu }_{-}e^+$ identify with each other's adjoints
over $\crnp$, because of the support preserving properties; more precisely,
$$
\Xi ^{\mu }_{+ }\colon \dot H_{p'}^{\operatorname{Re}\mu -s}(\crnp)\to \dot
 H_{p'}^{-s}(\crnp)  \text{ and }r^+\Xi ^{\overline\mu }_{-}e^+\colon \ol H_p^s(\rnp)\to
\ol H_p^{s-\operatorname{Re}\mu }(\rnp)\text{ are adjoints},
\tag A.8
$$
for $1<p<\infty $ and $1/p+1/p'=1$, all $s\in{\Bbb R}$. The same holds
for the operators $\Lambda _+^\mu ,\Lambda ^{\overline\mu }_{-,+}$, and there is a
similar statement for $\Lambda ^{(\mu )}_+$ and $\Lambda ^{(\overline
\mu )}_{-,+}$ relative to the set $\Omega $.

The following special spaces 
introduced by
H\"ormander \cite{H65} (for $p=2$), cf.\ \cite{G15}, are particularly adapted to $\mu
$-transmission operators $P$: 
$$
\aligned
\Cal E_\mu (\comega)&=e^+\{u(x)=d(x)^\mu v(x)\mid v\in C^\infty
(\comega)\},\\
H^{\mu (s)}_p(\crnp)&=\Xi _+^{-\mu }e^+\ol H_p^{s-\operatorname{Re}\mu
}(\rnp),\quad  s>\operatorname{Re}\mu -1/p',\\
H^{\mu (s)}_p(\comega)&=\Lambda  _+^{(-\mu )}e^+\ol H_p^{s-\operatorname{Re}\mu
}(\Omega ),\quad  s>\operatorname{Re}\mu -1/p'.
\endaligned\tag A.9
$$
Namely, $r^+P$ (of order $m$) maps them into $C^\infty (\comega)$,  $\ol H_p^{s-\operatorname{Re}m}(\rnp)$, resp.\ $\ol
 H_p^{s-\operatorname{Re}m}(\Omega)$ 
(cf.\
 \cite{G15} Sections 1.3, 2, 4).
In the first line of (A.9), $\operatorname{Re}\mu
 >-1$ (for other $\mu $, cf.\ \cite{G15}) and $d(x)$ is a $C^\infty $-function vanishing to
 order 1 at $\partial\Omega $ and positive on $\Omega $, e.g.\
 $d(x)=\operatorname{dist}(x,\partial\Omega )$ near $\partial\Omega $.

If in addition $P$ is elliptic with factorization index $\mu _0$
($\equiv \mu$ mod 1),
the Dirichlet problem
$$
r^+Pu=f,\quad \operatorname{supp}u\subset \comega, \tag A.10
$$
satisfies by \cite{G15}, Th.\ 4.4: When $u\in \dot H^\sigma _p(\comega)$
 (with $\sigma >\operatorname{Re}\mu _0-1/p'$) solves (A.10) for some
 $f\in \ol H^{s-m}_p(\Omega )$ with
 $s>\operatorname{Re}\mu _0-1/p'$, then $u\in H^{\mu _0(s)}_p(\comega)$;
 moreover, $r^+P$ is Fredholm from  $ H^{\mu _0(s)}_p(\comega)$ to $\ol
 H^{s-m}_p(\Omega )$. This will be used in the
 present paper with $\mu =\mu _0=a$, $m=2a$ for some $a\in \,]0,1[\,$.

 One has that $H^{\mu
 (s)}_p(\comega)\supset \dot H_p^s(\comega)$, and the distributions are
 locally in $H^s_p$   on $\Omega $, but at the boundary they in general have a
 singular behavior (cf.\ \cite{G15} Th.\ 5.4):
$$
H_p^{\mu (s)}(\comega)\cases =\dot H_p^s(\comega) \text{ if }s\in
\,]\operatorname{Re\mu }-1/p',\operatorname{Re}\mu +1/p[\,,\\
\subset e^+d^\mu  \ol H_p^{s- \operatorname{Re}\mu }(\Omega)+\dot
H_p^{s}(\comega)\text{ if }s>\operatorname{Re}\mu +1/p.
\endcases
\tag A.11
$$

\subhead Acknowledgement\endsubhead
We thank X.\ Ros-Oton for
useful conversations.

\enddocument